\documentclass{article}
\usepackage{graphicx}
\usepackage{authblk}
\usepackage{cite}

\usepackage{tikz}
\usetikzlibrary{arrows,calc}

\usepackage{graphicx}
\usepackage{float}
\usepackage{nicematrix}

\usepackage[colorlinks]{hyperref}

\usepackage{xcolor}

\usepackage{amsmath}
\usepackage{amssymb}
\usepackage{amsthm}

\usepackage{cleveref}

\newtheorem{theorem}{Theorem}
\newtheorem{lemma}[theorem]{Lemma}
\newtheorem{corollary}[theorem]{Corollary}
\newtheorem{proposition}[theorem]{Proposition}

\newtheorem{problem}{Problem}

\tikzstyle{arc}=[->,shorten <=3pt, shorten >=3pt, >=stealth, line width=1.1pt]
\tikzstyle{edge}=[shorten <=2pt, shorten >=2pt, >=stealth, line width=1.1pt]
\tikzstyle{blackV}=[circle, fill=black, draw, minimum size=5pt, inner sep=0pt,
                    outer sep=0pt]
\tikzstyle{vertex}=[circle, fill=white, draw, minimum size=5pt, inner sep=0pt,
                    outer sep=0pt]

\title{Strong chordality in tournaments and multipartite tournaments with
possible loops\thanks{The authors gratefully acknowledge support from NSERC
Canada, SEP-CONACYT CB-2017 A1-S-8397 grant and UNAM-PAPIIT IN106425 grant}}

\author[1]{Pavol~Hell\thanks{pavol@sfu.ca}}
\author[2]{C\'esar~Hern\'andez-Cruz\thanks{chc@ciencias.unam.mx}}
\author[3]{Jing~Huang\thanks{huangj@uvic.ca}}

\affil[1]{School of Computing Science,
          Simon Fraser University}
\affil[2]{Facultad de Ciencias,
          Universidad Nacional Aut\'onoma de M\'exico}
\affil[3]{Department of Mathematics and Statistics,
		      University of Victoria}

\begin{document}
\date{\today}

\maketitle
\begin{abstract}
    Strongly chordal digraphs are included in the class of chordal digraphs and
    generalize strongly chordal graphs and chordal bipartite graphs. They are
    the digraphs that admit a linear ordering of its vertex set for which their
    adjacency matrix does not contain the $\Gamma$ matrix as a submatrix. In
    general, it is not clear if these digraphs can be recognized in polynomial
    time. We focus on multipartite tournaments with possible loops. We give a
    polynomial-time recognition algorithm and a forbidden induced subgraph
    characterization of the strong chordality for each of the following cases:
    \begin{itemize}
        \item tournaments with possible loops,
        \item reflexive multipartite tournaments,
        \item irreflexive bipartite tournaments,
        \item irreflexive tournaments minus one arc, and
        \item balanced digraphs.
    \end{itemize}
    In addition, we prove that in a strongly chordal digraph the minimum size of
    a total dominating set equals the maximum number of disjoint
    in-neighborhoods, and this number can be calculated in linear time given a
    $\Gamma$-free ordering of the input graph.
\end{abstract}

\section{Background and definitions}
\label{sec:introduction}

A number of interesting graph classes have been extended to digraphs, including
interval graphs \cite{federDAM160}, chordal graphs \cite{haskinsSIAMJC2,
meisterTCS463,hy,hellDAM216}, and split graphs \cite{hellDAM216,lamarDM312}. In
most cases, there is more than one way to define such a generalization, and it
is not obvious which one best captures the analogy to the undirected case. (In
the undirected case there may be several equivalent characterizations of the
graphs in the class, and each may suggest a different generalization, which are
not equivalent in the context of digraphs.) It seems to be the case that often
the most successful generalizations use the ordering characterization of the
undirected concept, or, equivalently, its characterization by forbidden
submatrices of the adjacency matrix.

Consider first the undirected notion of an interval graph. Since every interval
intersects itself, we will assume each vertex has a loop. Then interval graphs
are known to have the following ordering characterization \cite{federDAM160}.
(There are other ordering characterizations of interval graphs, but this one
turns out to be most useful; however, it only applies if every vertex is
considered adjacent to itself.) A graph $G$ is an interval graph if and only if
its vertices can be ordered as $v_1, v_2, \dots, v_n$ so that if $i < j$ and $k
< \ell$, not necessarily all distinct, then for $v_iv_{\ell} \in E(G), v_jv_k
\in E(G)$ we also have $v_jv_{\ell} \in E(G)$. Equivalently, $G$ is an interval
graph if and only if the rows and columns of its adjacency matrix can be
simultaneously permuted to avoid a submatrix of the form
$\left[ \begin{smallmatrix}
  \ast & 1 \\
  1 & 0
\end{smallmatrix} \right]$
where $\ast$ can be either $0$ or $1$.

In \cite{federDAM160}, the authors analogously define a digraph analogue of
interval graphs as follows. A digraph with a loop at every vertex is an {\em
adjusted interval digraph} if the rows and columns of its adjacency matrix can
be simultaneously permuted to avoid a submatrix of the form
$\left[ \begin{smallmatrix}
  \ast & 1 \\
  1 & 0
\end{smallmatrix} \right].$

It turns out that these digraphs have a natural geometric representation, a
forbidden structure characterization, and other desirable properties analogous
to interval graphs \cite{federDAM160}. (By contrast, the earlier class of {\em
interval digraphs} \cite{dasJGT13}, based on a simple geometric analogy, lacks
many of these nice properties.)

To simplify the language, we will say that a vertex is {\em reflexive} if it has
a loop and {\em irreflexive} if it does not. A digraph is {\em reflexive} if
every vertex is reflexive and is {\em irreflexive} if every vertex is
irreflexive. Thus the diagonal entries of the adjacency matrix of a reflexive
digraph are all $1$ and of an irreflexive digraph are all $0$. An arc $(u,v)$ in
a digraph is {\em symmetric} if $(v,u)$ is also an arc. A digraph is {\em
symmetric} if every arc is symmetric. The adjacency matrix of a symmetric
digraph is symmetric. A symmetric digraph may be viewed as a graph with possible
loops. In the figures, we will depict reflexive vertices in {\em black} and
irreflexive vertices in {\em white}.

For graph classes that are characterized as intersection graphs (typically
chordal graphs and their subclasses such as strongly chordal graphs and interval
graphs), it is most natural to restrict attention to reflexive graphs (and
digraphs), as is noted above for interval graphs. Nevertheless, it is possible
to obtain useful generalizations for digraphs that are neither reflexive nor
irreflexive. This is done, for example, in \cite{arash,ross}, where general
digraphs (that have some vertices with loops and others without) avoiding
$\left[ \begin{smallmatrix}
  * & 1 \\
  1 & 0
\end{smallmatrix} \right]$
are investigated and found a useful unification of interval graphs, adjusted
interval digraphs, two-dimensional orthogonal ray graphs (alias interval
containment digraphs), and complements of threshold tolerance graphs.

In this paper we consider the digraph generalization of the undirected notion of
strong chordality. A chordal graph $G$ can be defined by the existence of a {\em
perfect elimination ordering}, also known as a {\em simplicial ordering}, $v_1,
v_2, \dots, v_n$ of its vertices so that if $i < j, i < k$ and $v_iv_j \in E(G),
v_iv_k \in E(G)$, then we must also have $v_jv_k \in E(G)$. They are also
characterized as those graphs that have no induced cycle of length greater than
three, or those graphs that are intersection graphs of subtrees of a tree
\cite{golumbic}. As noted above, we consider chordal graphs to be reflexive,
i.e., the adjacency matrix of a chordal graph has $1$'s on its main diagonal.
Then a perfect elimination ordering corresponds to a simultaneous permutation of
the rows and columns of the adjacency matrix that avoids as a principal
submatrix the so-called $\Gamma$ {\em matrix}
$\left[ \begin{smallmatrix}
  1 & 1 \\
  1 & 0
\end{smallmatrix} \right]$.
Such a submatrix is called a {\em principal submatrix} if the upper left $1$
lies on the main diagonal.

Chordal digraphs were first defined in \cite{haskinsSIAMJC2}, and further
studied in \cite{meisterTCS463}. A reflexive digraph $D$ is a {\em chordal
digraph} if the rows and columns of its adjacency matrix can be simultaneously
permutated to avoid $\Gamma$ as a principal submatrix. These digraphs can be
recognized in polynomial time \cite{haskinsSIAMJC2} and structural
characterizations are known for several special cases, including oriented graphs
and semi-complete digraphs \cite{meisterTCS463}, locally semicomplete digraphs
and weakly quasi-transitive digraphs \cite{hy}. A more restrictive notion of
{\em strict chordal digraphs} from \cite{hellDAM216} admits a general forbidden
induced subgraph characterization and leads to a nice notion of strict split
digraphs \cite{hellDAM216}.

In the context of undirected graphs, {\em strongly chordal graphs}
\cite{farberDM43} are defined as the subclass of those chordal graphs for which
the rows and columns of their adjacency matrix can be simultaneously permutated
to avoid $\Gamma$ as {\em any} submatrix (not just principal submatrix).
Strongly chordal graphs admit elegant forbidden structure characterizations
\cite{farberDM43,cn}, efficient recognition algorithms \cite{lubiwSIAMJC16}, and
lead to efficient algorithms for some problems that are intractable for chordal
graphs \cite{farberDM43}.

Permuting rows and columns of a $0, 1$ matrix $M$ to avoid $\Gamma$ as a
submatrix has been much studied \cite{ansteeJA5,hoffmanSIAMJADM6,lubiwSIAMJC16}.
A {\em $\Gamma$-free ordering} of $M$ is a matrix obtained from $M$ by {\em
independently} permuting its rows and columns, to avoid $\Gamma$ as a submatrix.
If the constraint matrix of a linear program is presented in a $\Gamma$-free
ordering, then it can be solved by a greedy algorithm
\cite{ansteeJA5,hoffmanSIAMJADM6}. A {\em cycle matrix} is a square $0, 1$
matrix of size at least $3$, with exactly two $1$'s in each row and each column.
A matrix $M$ is {\em totally balanced}, if it admits no cycle matrix as a
submatrix.  A matrix $M$ admits a $\Gamma$-free ordering if and only if it is
totally balanced \cite{hoffmanSIAMJADM6}.  There are efficient algorithms to
decide if a matrix is totally balanced \cite{lubiwSIAMJC16}.

For a square matrix $M$, a {\em symmetric $\Gamma$-free ordering} is a matrix
obtained from $M$ by {\em simultaneously} permuting its rows and columns, to
avoid $\Gamma$ as a submatrix. A reflexive graph $G$ is strongly chordal if and
only if its adjacency matrix $M(G)$ has a symmetric $\Gamma$-free ordering
\cite{farberDM43}. The algorithm in \cite{lubiwSIAMJC16} finds a symmetric
$\Gamma$-free ordering of a symmetric matrix $M$ (or decides that one doesn't
exist) provided $M$ has $1$'s on the main diagonal. In particular, a symmetric
matrix $M$ with $1$'s on the main diagonal admits a symmetric $\Gamma$-free
ordering if and only if it is totally balanced \cite{farberDM43,lubiwSIAMJC16}.

For a bigraph $G$ (a bipartite graph with a fixed bipartition into red and blue
vertices), we consider the {\em bi-adjacency matrix} $N(G)$, with rows indexed
by the red vertices and columns indexed by the blue vertices, and $N(i, j)=1$ if
and only if the $i$-th red vertex is adjacent to the $j$-th blue vertex. Note
that $N$ is in general not a square matrix. A {\em chordal bigraph} $G$ is a
bigraph whose bi-adjacency matrix has a $\Gamma$-free ordering \cite{hmp}.

We say $D$ is {\em a strongly chordal digraph} if its adjacency matrix $M(D)$
admits a symmetric $\Gamma$-free ordering. It follows that a strongly chordal
graph is precisely (the underlying graph of) a strongly chordal digraph that is
symmetric and reflexive. It also follows that strongly chordal digraphs are
chordal digraphs as defined in \cite{haskinsSIAMJC2,meisterTCS463}. Chordal
bigraphs can also be seen as special strongly chordal digraphs, because the
adjacency matrix $M(G)$ of a bigraph $G$ (viewed as a graph) has a symmetric
$\Gamma$-free ordering if and only if its bi-adjacency matrix $N(G)$ has a
$\Gamma$-free ordering. Thus strongly chordal digraphs can be seen as
generalizing strongly chordal graphs, and chordal bigraphs, and be included in
the class of chordal digraphs. The strong chordality of symmetric digraphs with
possible loops has been characterized in \cite{hhhl}.

The problem of recognizing strongly chordal digraphs is equivalent to the
problem of deciding if a given square $0, 1$ matrix has a symmetric
$\Gamma$-free ordering. This seems to be a difficult problem; as we show below,
it is no longer equivalent with being totally balanced, or any of the other
polynomial conditions that applied for symmetric matrices with $1$'s on the main
diagonal.

We shall focus on certain particular classes of digraphs. The first class is the
class of tournaments with possible loops in \Cref{sec:tournaments}; here we
prove that very few of these tournaments are strongly chordal, and we can
actually describe all strongly chordal cases.   Irreflexive bipartite
tournaments are considered in \Cref{sec:irr-bip}, where we obtain
characterizations for strongly chordal members of this class in terms of
forbidden induced subdigraphs.   We use irreflexive tournaments minus one arc in
\Cref{sec:irr-minus} as an example to show that the problem of characterizing
strongly chordal irreflexive multipartite tournaments is already a challenging
problem; we also obtain a characterization of the strongly connected members of
this family in terms of forbidden induced subdigraphs.   We exhibit
characterizations of strongly chordal reflexive multipartite tournaments in
\Cref{sec:ref-mult}, both by describing their global structure and in terms of
forbidden induced subdigraphs.   Results in
\Cref{sec:tournaments,sec:irr-bip,sec:irr-minus,sec:ref-mult} are used to
provide recognition algorithms for the discussed families running in
$O(|V|+|A|)$ which are certifying, except for the case of irreflexive bipartite
tournaments.   In \Cref{sec:application} we discuss a concept that generalizes
both domination and total domination in digraphs when loops and non-loops are
allowed, and show that the associated parameter can be computed in polynomial
time given a $\Gamma$-free ordering of a strongly chordal digraph. Conclusions
and open problems are addressed in \Cref{sec:conclusions}; in this final section
we also consider strongly chordal balanced digraphs, which are a different
generalization of chordal bigraphs and include all oriented trees. But first, we
introduce nomenclature and prove a number of simple yet useful results in
\Cref{sec:general}.

\section{General Digraphs}
\label{sec:general}

We consider digraphs where loops are allowed but without parallel arcs.   A
digraph $D$ is considered to have vertex set $V(D)$ and arc set $A(D)$.  For a
pair of vertices $u$ and $v$ of a digraph $D$, we sometimes write $u \to v$ to
denote that $(u,v)$ is an arc of $D$.

A vertex $v$ in a digraph $D$ is {\em simplicial} if it is either irreflexive or
reflexive and for all vertices $u \in N^-(v)$ and $w \in N^+(v)$, there is an
arc $(u,w) \in A(D)$. A {\em simplicial ordering} of a digraph $D$ is a linear
ordering $v_1, v_2, \dots, v_n$ of its vertices, such that for each $i$, the
vertex $v_i$ is simplicial in $D - \{v_1, v_2, \dots, v_{i-1}\}$. A digraph is
{\em chordal} if and only if it has a simplicial ordering.

We will call a vertex $v$ in a digraph $D$ {\em simple} if
\begin{itemize}
  \item $v$ is simplicial,
  
  \item if $x, y \in N^-(v)$, then $N^+(x) \subseteq N^+(y)$ or $N^+(x)
        \supseteq N^+(y)$, and

  \item if $x, y \in N^+(v)$, then $N^-(x) \subseteq N^-(y)$ or $N^-(x)
        \supseteq N^-(y)$.
\end{itemize}

A {\em simple ordering} of a digraph $D$ is a vertex ordering $v_1, v_2, \dots,
v_n$ of $D$ such that for each $i$, the vertex $v_i$ is simple in $D - \{v_1,
v_2, \dots, v_{i-1}\}$. Observe that a simple ordering is again a simplicial
ordering. A {\em strong ordering} of a digraph $D$ is a linear ordering $v_1,
v_2, \dots, v_n$ of its vertices such that for all $i < j$ and $k < \ell$ where
$i, j, k, \ell$ are not necessarily distinct, if $(v_i,v_k), (v_i,v_{\ell}),
(v_j,v_k) \in A(D)$, then $(v_j,v_{\ell}) \in A(D)$. A strong ordering of $D$
directly corresponds to a symmetric $\Gamma$-free ordering of $M(D)$. A strong
ordering is a simple ordering and hence a simplicial ordering. A digraph is
strongly chordal if and only if it has a strong ordering. Thus each strongly
chordal digraph is a chordal digraph.

Having a simple ordering is equivalent to having a strong ordering in the
classical context (reflexive symmetric digraphs), but is not equivalent for
general digraphs. Every symmetric $\Gamma$-free ordering is a simple ordering,
but the converse is not necessarily true. This is not true even for irreflexive
tournaments; the irreflexive tournament $T_1$ in Figure~\ref{fig:minobsirrtour}
has a simple ordering (and its adjacency matrix is totally balanced), but is not
strongly chordal, i.e., the matrix has no symmetric $\Gamma$-free ordering.  The
relation between simple and strong orderings yields a simple necessary condition
for a digraph to be strongly chordal.

\begin{lemma}
\label{lem:simple}
    Let $D$ be a digraph.   If no vertex of $D$ is simple, then $D$ is not
    strongly chordal.
\end{lemma}

\begin{proof}
    As we noted in the previous paragraph, every strongly chordal graph has a
    simple vertex, which is actually the first vertex of the strong ordering.
    Hence, if no vertex of $D$ is simple, then no vertex of $D$ can be the first
    in the strong ordering, and thus it is not strongly chordal.
\end{proof}

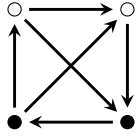
\begin{figure}[htb!]
\begin{center}
\begin{tikzpicture}
    \node [blackV] (1) at (0,0) {};
    \node [vertex] (2) at (0,1.5) {};
    \node [vertex] (3) at (1.5,1.5) {};
    \node [blackV] (4) at (1.5,0) {};
    \foreach \from/\to in {1/2,1/3,2/3,2/4,3/4,4/1}
        \draw [arc] (\from) to (\to);
\end{tikzpicture}
\end{center}
\vspace{-2mm}
\caption{Tournament $T_0$.}
\label{mixed_tournament}
\end{figure}

The tournament $T_0$ in \Cref{mixed_tournament} contains both reflexive and
irreflexive vertices. It follows from \Cref{lem:simple} that $T_0$ is not
strongly chordal although each of the subgraphs induced by reflexive and
irreflexive vertices respectively is strongly chordal.

The necessary condition given in \Cref{lem:simple} is often not enough to verify
that a digraph is not strongly chordal.   We present a series of results
describing properties that every strong ordering of a digraph must fulfill.

\begin{lemma}
\label{lem:ref-or-P3}
    Let $D$ be a reflexive oriented graph, and let $u, v$ and $w$ be vertices of
    $D$. If $\{ u, v, w \}$ induces an orientation of $P_3$ having $v$ as its
    middle vertex, then $v$ cannot appear before $u$ and $w$ in an strong
    ordering of $D$.
\end{lemma}

\begin{proof}
    If $(u,v,w)$ is a directed path, then $v$ is not a simplicial vertex.  Since
    $u$ is an in-neighbour and an out-neighbour of itself and it is neither an
    in-neighbour nor an out-neighbour of $w$, and vice versa, if  $u \to v \gets
    w$ or $u \gets v \to w$, then $v$ is not a simple vertex.   Our previous
    observations remain true in any induced subdigraph of $D$ where $u, v$ and
    $w$ are part of the vertex set, therefore, $v$ cannot be the first vertex in
    a strong ordering of any such subdigraph.
\end{proof}

\begin{lemma}
\label{lem:simple-ord}
    Let $D$ be an oriented graph with a $\Gamma$-free ordering $\le$, and let
    $x$ be a vertex of $D$.
    \begin{enumerate}
        \item If $y,z \in N^+(x)$ appear after $x$ in a strong ordering of $D$
            and $N_{D_x}^-(y) \subset N_{D_x}^-(z)$, then $y$ is placed before
            $z$ in the strong ordering.
        \item If $y,z \in N^-(x)$ appear after $x$ in a strong ordering of $D$
            and $N_{D_x}^+(y) \subset N_{D_x}^+(z)$, then $y$ is placed before
            $z$ in the strong ordering.
    \end{enumerate}
\end{lemma}

\begin{proof}
    Due to their similarity, we only prove the first item.   Clearly, $x$ is the
    first vertex in the $\Gamma$-free ordering of $D_x$ resulting from
    restricting the $\Gamma$-free ordering of $D$ to $D_x$, so it is simple in
    $D_x$. Hence, for $y,z \in N_{D_x}^+(x)$, we may assume without loss of
    generality that $N_{D_x}^-(y) \subseteq N_{D_x}^-(z)$.   If $N_{D_x}^-(y)
    \subset N_{D_x}^-(z)$, then there is $w \in N_{D_x}^-(z) \setminus
    N_{D_x}^-(y)$, so having $z \le y$ would create a $\Gamma$ in the submatrix
    of $A_D$ induced by the rows corresponding to $x$ and $w$, and the columns
    corresponding to $y$ and $z$ (see \Cref{fig:simple-ord}).
    \begin{figure}[ht!]
        \centering
        \[
            \begin{pNiceArray}{ccccc}[first-col,first-row]
                  &        & z      &        & y      & \\
                  &        & \Vdots &        & \Vdots & \\
                x & \Cdots & 1      & \Cdots & 1      & \\
                  &        & \Vdots &        & \Vdots & \\
                w & \Cdots & 1      & \Cdots & 0      & \\
                  &        &        &        &        &
            \end{pNiceArray}
        \]
    \caption{A $\Gamma$ submatrix created by an incorrect order of the later
    out-neighbours of a simple vertex.}
    \label{fig:simple-ord}
    \end{figure}
\end{proof}

A vertex $v$ in a digraph $D$ is a {\em peak} vertex of $D$ if there exist
vertices $u, w \in V(D)$ such that $(u,v) \in A(D), (v,w) \in A(D)$ and $(u,w)
\in A(D)$. For irreflexive digraphs, we can forbid some vertices to be the last
in a strong ordering.

\begin{lemma}\label{lem:peak}
An irreflexive vertex that is a peak cannot be the last vertex in a simple
ordering.
\end{lemma}

\begin{proof}
Let $v_1, \dots, v_n$ be a simple ordering of the vertices of $D$, and assume
$v_n$ is irreflexive and a peak vertex with arcs $v_iv_j, v_iv_n, v_nv_j$ in
$G$. Then a $\Gamma$ submatrix occurs in rows $i, n$ and columns $j, n$.
\end{proof}

\begin{corollary}
Let $D$ be an irreflexive digraph. If every vertex of $D$ is a peak, then $D$ is
not strongly chordal.
\end{corollary}

If follows from the proofs of Lemmas \ref{lem:simple} and \ref{lem:peak} that in
any strong ordering of a strongly chordal digraph $D$ the first vertex must be
simple and the last vertex must not be an irreflexive peak. Therefore, we
observe for future reference that if an irreflexive digraph $D$ has only one
vertex that is simple, and at the same time it is the only vertex of $D$ which
is not a peak, then $D$ is not strongly chordal.

We now turn our attention to the irreflexive case, where we present some simple
conditions making an oriented graph not strongly chordal.   The following
technical lemma is quite simple, but will be useful to deal with some
cornerstone cases in our analysis.

\begin{lemma}
\label{lem:irr-two-arcs}
    Let $D$ be an irreflexive oriented graph, let $(u,v)$ and $(x,y)$ be two
    different arcs of $D$ such that $u \not\to y$ and $x \not\to v$, and let $w$
    be a vertex of $D$.   If $w$ dominates $v$ and $y$, or $u$ and $x$ dominate
    $w$, then $w$ is not simple.
\end{lemma}
\begin{proof}
    It follows from our hypotheses that $u \ne x$ and $v \ne y$, but it is
    possible that $v = x$.   Suppose first that $v = x$.   In this case, $x$ and
    $y$ are out-neighbours of $w$ such that $x$ is an in-neighbour of $y$ but
    not of itself and $u$ is an in-neighbour of $x$ but not of $y$.   Thus, $x$
    and $y$ have incomparable out-neighbourhoods, and we conclude that $w$ is
    not simple.   Suppose now that $v \ne x$.   By hypothesis, $u$ is an
    in-neighbour of $v$ which is not an in-neighbour of $y$, and $x$ is an
    in-neighbour of $y$ which is not an in-neighbour of $v$.  Hence, $w$ has two
    out-neighbours with incomparable in-neighbourhoods, so it is not simple.
    The remaining case has a similar proof.
\end{proof}

\begin{lemma}
\label{lem:irr-cyc}
    Let $D$ be an irreflexive oriented graph, and let $v$ be a vertex of $D$.
    If any of the following conditions hold:
    \begin{enumerate}
        \item $v$ has two in-neighbours or two out-neighbours in an induced
            directed cycle of $D$,
        \item $v$ has two out-neighbours different from the initial vertex in
            an induced directed path of $D$,
        \item $v$ has two in-neighbours different from the terminal vertex in an
            induced directed path of $D$,
    \end{enumerate}
    then $v$ is not simple.
\end{lemma}
\begin{proof}
    All cases follow directly from \Cref{lem:irr-two-arcs} by considering the
    two arcs defined by the two out-neighbours of $v$ and their in-neighbours in
    the path or cycle, or the two arcs  defined by the two in-neighbours of $v$
    and their out-neighbours in the path or cycle.
\end{proof}

\begin{proposition}
    Let $H$ be an irreflexive oriented graph, and let $D$ be obtained from $H$
    by adding a new vertex $x$ dominating every vertex of $H$, and a new vertex
    $y$ dominated by every vertex of $H$.   The arc from $x$ to $y$ is missing
    in $D$, but the arc from $y$ to $x$ might be present.   If $H$ is any of the
    following:
    \begin{enumerate}
        \item a directed $3$-cycle,
        \item a directed $P_3$,
        \item an oriented $2K_2$,
    \end{enumerate}
    then $D$ is a minimal obstruction for strong chordality.
\end{proposition}
\begin{proof}
    For the first two items it follows from \Cref{lem:irr-cyc} that $x$ and $y$
    are not simple, and for the last item we reach the same conclusion by
    applying \Cref{lem:irr-two-arcs}.   Now, let $u$ be any vertex in $H$.   If
    $u$ has two out-neighbours, then one of them, say $v$, is also in $H$, and
    the other one is $y$.   Since the arc from $x$ to $y$ is missing by
    hypothesis, and $u$ dominates $v$ and $y$, we apply \Cref{lem:irr-two-arcs}
    to $u$ and the arcs $(x,v)$ and $(v,y)$ to conclude that $u$ is not simple.
    Otherwise, $u$ has two in-neighbours.  Again, one of them, say $v$ is in
    $H$, and the other one is $x$.   Analogous to the previous case, we apply
    \Cref{lem:irr-two-arcs} to $u$ and the arcs $(x,v)$ and $(v,y)$ to conclude
    that $u$ is not simple.   Therefore, no vertex of $D$ is simple, and
    \Cref{lem:simple} implies that $D$ is not strongly chordal.

    The minimality for the first two items is easy to establish, as deleting any
    vertex yields a tournament on $4$ vertices or the transitive tournament on
    $4$ vertices minus the arc from the source to the sink.   For the last item,
    it is clear from the fact that the adjacency matrix of $H$ has only two
    ones, that when the source or the sink is deleted any ordering it suffices
    to put the remaining sink or source at the end to obtain a $\Gamma$-free
    ordering.  If a sink $v$ in $H$ is deleted from $D$, then a $\Gamma$-free
    ordering of $D-v$ is obtained by starting with the remaining isolated vertex
    in $H$, following with the sink of $D$, the sink and source (in that order)
    of the remaining arc in $H$, and finishing with the source of $D$.   A
    $\Gamma$-free ordering for the case when a source in $H$ is deleted from $D$
    is similarly constructed.
\end{proof}

We finish this section with a result relating simple orderings with totally
balanced matrices.

\begin{proposition}\label{jing}
If a digraph $D$ has a simple ordering, then $M(D)$ is totally balanced.
\end{proposition}

\begin{proof}
Suppose a digraph $D$ has a simple ordering. A digraph $D$ defines a bigraph
$B(D)$ (just as for graphs): each vertex $v \in V(D)$ gives rise to two vertices
$v_1, v_2 \in V(B(D))$, and each arc $(v,w) \in A(D)$ gives rise to an edge
$(v_1,w_2)$ of $B(D)$. Then it is easy to see that the bigraph $B(D)$ also has a
simple ordering. Thus the bi-adjacency matrix of $B(D)$ is totally balanced.
Moreover, we have $N(B(D))=M(D)$.
\end{proof}

\section{Tournaments}
\label{sec:tournaments}

As we have seen, strongly chordal digraphs do not in general coincide with
digraphs having a simple ordering, or having a totally balanced adjacency
matrix, even for tournaments.

We begin by addressing two natural subcases, reflexive and irreflexive
tournaments.   Clearly, the matrix of a reflexive directed cycle on three
vertices is not totally balanced (it is itself the bi-adjacency matrix of an
even cycle of length $6$). Thus, every reflexive strongly chordal tournament is
acyclic, and we have the following theorem.

\begin{theorem}
\label{thm:reftour}
If $T$ is a reflexive tournament, then $T$ is strongly chordal if and only if it
is a reflexive transitive tournament.
\end{theorem}

The irreflexive case, although more interesting, is similar in flavor to the
reflexive one.   Let us start with a simple yet useful observation.

\begin{lemma}
\label{lem:irr-sin-sou}
    Let $T$ be an irreflexive strongly chordal tournament.   If $T'$ is a
    tournament obtained from $T$ by adding a source or a sink, then $T'$ is
    strongly chordal.
\end{lemma}
\begin{proof}
    In order to obtain a $\Gamma$-free ordering for the vertices of $T'$, it
    suffices to consider a $\Gamma$-free ordering for the vertices of $T$, and
    add the new vertex at the end.  
\end{proof}

We now state an elementary necessary condition for an irreflexive tournament to
be strongly chordal.

\begin{lemma}
\label{lem:irr-two-cyc}
    Let $T$ be an irreflexive tournament.   If $T$ contains two directed
    triangles not sharing an arc, then $T$ is not strongly chordal.
\end{lemma}
\begin{proof}
    Let $C_1$ and $C_2$ be the two directed triangles.   Suppose first that $C_1
    = (u,v,w,u)$ and $C_2 = (w,x,y,w)$.   Observe that both $u$ and $v$ have at
    least two in-neighbours or at least two out-neighbours in $C_2$, and
    similarly, $x$ and $y$ have at least two in-neighbours or at least two
    out-neighbours in $C_1$.   Hence, following \Cref{lem:irr-cyc}, none of $u,
    v, x, y$ is simple.
    
    Since $T$ is a tournament, if $\{u, v, x, y\}$ contains a directed triangle,
    then $w$ has at least two in-neighbours or at least two out-neighbours in
    such a cycle, and \Cref{lem:irr-cyc} implies that it is not simple.   In
    this case, $\{u, v, w, x, y\}$ induces a subdigraph of $T$ without simple
    vertices, \Cref{lem:simple} implies that $T$ is not strongly chordal.
    
    So, suppose that $\{u, v, x, y\}$ induces a transitive tournament with
    topological ordering $(u, v, x, y)$.   Let $T'$ be the subtournament of $T$
    induced by $\{u, v, w, x, y\}$.   Clearly, $v$ and $x$ are peak vertices in
    $T'$.   Notice that $w$ dominates $u$ and is dominated by $y$, as otherwise
    the former would be a source in $T'$, and the latter a sink in $T'$, which
    is impossible because both of them belong to directed triangles in $T'$. But
    now, $w$ dominates $u$ and either $v$ or $x$, which makes $u$ a peak vertex,
    and $w$ is also dominated by $y$ and either $v$ or $x$, which makes $y$ a
    peak vertex.   Therefore, $w$ is the only simple and non-peak vertex in
    $T'$, and hence, $T$ is not strongly chordal.

    Suppose now that $C_1 = (u,v,w,u)$ and $C_2 = (x,y,z,x)$ are vertex
    disjoint.   If every vertex of $C_1$ dominates every vertex of $C_2$, or
    vice versa, then it follows from \Cref{lem:irr-cyc} that no vertex in the
    subtournament induced by $\{u, v, w, x, y, z\}$ is simple.   Thus, we deduct
    from \Cref{lem:simple} that $T$ is not strongly chordal.   Else, suppose
    without loss of generality that $w$ has an in-neighbour and an out-neighbour
    in $C_2$.   It is easy to verify that there are two vertices in $C_2$ which
    together with $w$ induce a directed triangle.   Hence, by considering such
    vertices together with $u, v$ and $w$, we are back in the case where we have
    two arc-disjoint triangles sharing a vertex.
\end{proof}

A first application of \Cref{lem:irr-two-cyc} is to exhibit some tournament
minimal obstructions for strong chordality.   But first, let us introduce a
family of strongly connected and strongly chordal irreflexive tournaments.

For every integer $n$, $n \ge 3$, let $TT_n$ and $TT^\ast_n$ denote the
irreflexive transitive tournament on $n$ vertices, and the tournament obtained
from the irreflexive transitive tournament on $n$ vertices where the arc from
the only source to the only sink has been reversed.   It is easy to verify that
ordering the vertices of $TT_n$ increasingly with respect to their in-degrees
results in a $\Gamma$-free ordering; the same order, up to reversing the arc
from the first to the last vertices, is a $\Gamma$-free ordering for
$TT^\ast_n$.   Hence, $TT_n$ and $TT^\ast_n$ are strongly chordal digraphs for
every $n \ge 3$.

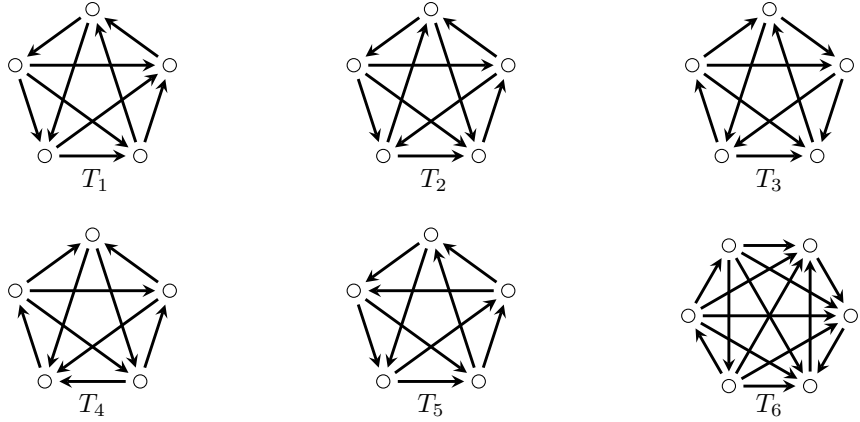
\begin{figure}[htb!]
\centering
\begin{tikzpicture}
\begin{scope}[scale=0.6]

\begin{scope}[xshift=-7.5cm,scale=0.9,rotate=162]
\foreach \i in {0,...,4}
	\draw ({(360/5)*\i}:2) node[vertex](\i){};
\foreach \i/\j in {0/1,0/2,0/3,1/2,1/3,2/3,2/4,3/4,4/0,4/1}
	\draw [arc] (\i) -- (\j);
\node (T1) at (110:2.2){$T_1$};
\end{scope}

\begin{scope}[xshift=0cm,scale=0.9,rotate=162]
\foreach \i in {0,...,4}
	\draw ({(360/5)*\i}:2) node[vertex](\i){};
\foreach \i/\j in {0/1,0/2,0/3,1/2,1/4,2/3,3/1,3/4,4/0,4/2}
	\draw [arc] (\i) -- (\j);
\node (T2) at (110:2.2){$T_2$};
\end{scope}

\begin{scope}[xshift=7.5cm,scale=0.9,rotate=234]
\foreach \i in {0,...,4}
  \draw ({(360/5)*\i}:2) node[vertex](\i){};
\foreach \i/\j in {0/1,0/4,1/3,2/0,2/1,3/0,3/2,4/1,4/2,4/3}
  \draw [arc] (\i) to (\j);
\node (T3) at (36:2.2){$T_3$};
\end{scope}

\begin{scope}[xshift=-7.5cm,yshift=-5cm,scale=0.9,rotate=306]
\foreach \i in {0,...,4}
	\draw ({(360/5)*\i}:2) node[vertex](\i){};
\foreach \i/\j in {0/1,0/4,1/2,1/4,2/0,2/4,3/0,3/1,3/2,4/3}
	\draw [arc] (\i) -- (\j);
\node (T4) at (324:2.2){$T_4$};
\end{scope}

\begin{scope}[xshift=0cm,yshift=-5cm,scale=0.9,rotate=90]
\foreach \i in {0,...,4}
	\draw ({(360/5)*\i}:2) node[vertex](\i){};
\foreach \i in {0,...,4}
	\draw [arc] let \n1={int(mod(\i+1,5))}
		in (\i) -- (\n1);
\foreach \i in {0,...,4}
	\draw [arc] let \n1={int(mod(\i+2,5))}
		in (\i) -- (\n1);
\node (T5) at (-2.2,0){$T_5$};
\end{scope}

\begin{scope}[xshift=7.5cm,yshift=-5cm,scale=0.9]
\foreach \i in {0,...,5}
  \draw ({(360/6)*\i}:2) node[vertex](\i){};
\foreach \i/\j in {0/5,1/0,2/0,2/1,2/4,2/5,3/0,3/1,3/2,3/5,4/0,4/1,4/3,4/5,5/1}
  \draw [arc] (\i) -- (\j);
\node (T6) at (0,-2.2){$T_6$};
\end{scope}

\end{scope}
\end{tikzpicture}
\caption{The family $\cal{T}$ of minimal obstructions to strong chordality.}
\label{fig:minobsirrtour}
\end{figure}

Let $\cal{T}$ be the family of tournaments $\{ T_1, \dots, T_6 \}$ depicted in
Figure \ref{fig:minobsirrtour}.

\begin{proposition}
\label{pro:irr-min-obs}
    Every tournament in the family $\mathcal{T}$ shown in
    \Cref{fig:minobsirrtour} is a minimal obstruction to strong chordality.
\end{proposition}
\begin{proof}
    It follows from \Cref{lem:irr-two-cyc} that every tournament in the family
    $\mathcal{T}$ is not strongly chordal.   There are exactly four tournaments
    of order $4$, namely $TT_4, TT_4^\ast$ and the tournaments obtained from a
    directed triangle by adding either a sink or a source.   We have already
    observed that the first two are strongly chordal, and the fact that the last
    two are strongly chordal follows from \Cref{lem:irr-sin-sou}.   Therefore,
    we have already verified that $T_i$ is a minimal obstruction for $i \in \{
    1, \dots, 5 \}$.

    For $T_6$, notice that the deletion of any vertex yields a tournament
    isomorphic to the tournament obtained by a directed triangle by adding
    either two sinks or two sources.   Again, \Cref{lem:irr-sin-sou} implies
    that the obtained tournament is strongly chordal, and hence, $T_6$ is a
    minimal obstruction to strong chordality.
\end{proof}

It is not hard to verify that the only strong tournaments on three and four
vertices are precisely $TT^\ast_3$ and $TT^\ast_4$, both of which are strongly
chordal.  Our next result extends this observation to tournaments of any order.

\begin{lemma}
\label{lem:irr-str-tour}
Let $n$ be an integer, $n \ge 3$.   The only irreflexive tournament on $n$
vertices which is both strongly connected and strongly chordal is $TT^\ast_n$.
\end{lemma}
\begin{proof}
By induction on $n$.   We have already noticed that the statement is true for $n
\in \{ 3, 4 \}$, so let $n$ be at least $5$, and let $T$ be a strongly connected
and strongly chordal tournament.   Since $T$ is strong, it is pancyclic, and
hence it contains a vertex $v$ such that the subtournament $T'$ obtained from
$T$ by deleting $v$ is strong.   Strong chordality is a hereditary property and
thus, by induction hypothesis, $T'$ is isomorphic to $TT^\ast_{n-1}$. Consider
an ordering $v_1, v_2, \dots, v_{n-2}, v_{n-1}$ of $V(T')$ such that $v_{n-1}
\to v_1$ and the reversal of this arc results in the transitive tournament
$TT_{n-1}$.

Since $T$ is strong, arcs from $T'$ to $v$ and from $v$ to $T'$ must exist in
$T$.   Let us observe that if $v \to v_i$, then $v \to v_j$  for every $i < j$.
Suppose that the contrary happens, this is, suppose that there are $i, j \in \{
1, \dots, n-1 \}$ such that $i < j$ and $v_j \to v \to v_i$.   Moreover, suppose
that $j-i$ is minimum among pairs $i,j$ with this property.   By the minimality
of $j-i$, either $i \ne 1$ or $j \ne n-1$; suppose without loss of generality
that $j \ne n-1$.   If $i \ne 1$, then $(v_1, v_i, v_{n-1}, v_1)$ and $(v_i,
v_j, v, v_i)$ are two arc-disjoint triangles in $T$.   If $i = 1$, then for any
$k \in \{ 2, \dots, n-2 \}$ different from $i$ and $j$ we have the arc-disjoint
triangles $(v_1,v_k,v_{n-1},v_1)$ and $(v_1, v_j, v, v_1)$.   In either case,
\Cref{lem:irr-two-cyc} implies that $T$ is not stronlgy chordal, a
contradiction.   Therefore, there exists $i \in \{ 2, \dots, n-2 \}$ such that
$\{ v_1, \dots, v_i \} \to v \to \{ v_{i+1}, \dots, v_{n-1} \}$.   Hence, it is
clear that reversing the arc $(v_{n-1}, v_1)$ in $T$ results in a transitive
tournament, and therefore $T$ is isomorphic to $TT^\ast_n$.
\end{proof}

Thus, in the strongly connected case, the only strongly chordal irreflexive
tournaments are very close to a transitive tournament. As the following argument
shows, this is also the case for general tournaments.

\begin{proposition}
\label{pro:irr-arc-rev}
    An irreflexive tournament is strongly chordal if and only if it is obtained
    from a transitive tournament by reversing at most one arc.
\end{proposition}
\begin{proof}
    We already know that transitive tournaments are strongly chordal, so it
    suffices to consider the case when $T$ is not transitive.   Suppose first
    that $T$ is a tournament obtained from the transitive tournament $T'$ by
    reversing exactly one arc.   So, let $v_1, \dots, v_n$ be a topological
    ordering of $T'$ and let $i, j \in \{ 1, \dots, n \}$ be such that $(v_i,
    v_j)$ is the arc reversed in $T'$ to obtain $T$. Since reversing an arc of
    the form $(v_i, v_{i+1})$ results in a transitive tournament, we may assume
    that $i+1 < j$.   Clearly, the tournament induced by $\{v_i, \dots, v_j\}$
    is isomorphic to $TT_{j-i}^\ast$, and $T$ is obtained from this tournament
    by adding sinks and sources.   Since $TT_{j-i}^\ast$ is strongly chordal, it
    follows from \Cref{lem:irr-sin-sou} that $T$ is strongly chordal.

    Conversely, let $T$ be a strongly chordal tournament.   If $T$ is strong,
    the result follows from \Cref{lem:irr-str-tour}.   Else, the fact that $T_6$
    is a minimal obstruction for strong chordality implies that $T$ contains at
    most one non-trivial strong component, say $C$.   Now, it follows from
    \Cref{lem:irr-str-tour} that $C$ is isomorphic to $TT_k^\ast$ for some $k$
    greater than or equal to $3$, and since every other connected component of
    $T$ is a single vertex, we conclude that $T$ is obtained from $C$ by adding
    some sinks and sources.   Clearly, $C$ has an arc such that its reversal
    causes $T$ to transform into a transitive tournament.
\end{proof}

In addition to the nice simple structure that irreflexive strongly chordal
tournaments have, it is possible to characterize them by a small set of minimal
forbidden induced subdigraphs.   Our next theorem provides two additional
characterizations for irreflexive strongly chordal tournaments, one in terms of
minimal forbidden induced subdigraphs, and we also show that the necessary
condition presented in \Cref{lem:irr-two-cyc} is also sufficient.

\begin{theorem}
\label{thm:irr-tou-cha}
    Let $T$ be an irreflexive tournament.   The following statements are
    equivalent. 
    \begin{enumerate}
        \item $T$ is strongly chordal.
        \item $T$ is obtained from a transitive tournament by reversing at most
            one arc.
        \item $T$ does not contain two arc-disjoint directed triangles.
        \item $T$ is $\cal{T}$-free.
    \end{enumerate}
\end{theorem}

\begin{proof}
Notice first that \Cref{pro:irr-arc-rev} proves the equivalence between the
first two items, \Cref{pro:irr-min-obs} proves that the first item implies the
fourth one, and \Cref{lem:irr-two-cyc} proves that the first item implies the
third one.

We prove that the last item implies the third one by contrapositive.   Suppose
that $T$ contains two arc-disjoint directed triangles, say $C_1$ and $C_2$.   If
the vertex sets of $C_1$ and $C_2$ are disjoint, then either $C_1 \to C_2$, and
hence $V_{C_1} \cup V_{C_2}$ induces $T_6$ or, as in the proof of
\Cref{lem:irr-two-cyc}, we find two arc-disjoint directed cycles sharing a
vertex.   So, from now on, we assume that $C_1 = (u, v, w, u)$ and $C_2 = (w, x,
y, w)$.   Let $T_0$ be the subtournament of $T$ induced by $\{u, v, x, y\}$, and
notice that $w$ has two out-neighbours and two in-neighbours in $T_0$.   There
are exactly four non-isomorphic tournaments on four vertices, and it is not hard
to notice that when $T_0$ is either $TT_4$, or it is obtained from $TT_3$ by
adding either a source or a sink, there is only one way in which a fifth vertex
can be added to create two arc-disjoint triangles (see $T_1, T_2, T_3$ in
\Cref{fig:minobsirrtour} where the fifth vertex is the top vertex).   When $T_0$
is the only strong tournament of order $4$, there are two ways of adding a fifth
vertex in such a way that two arc-disjoint triangles are created (see $T_4$ and
$T_5$ in \Cref{fig:minobsirrtour}), distinguished by the relative position of
the arcs of the new vertex to the diagonals of the only hamiltonian cycle in
$TT_4$.   Hence, we conclude that if $T$ is not $\mathcal{T}$-free, as desired.

Finally, let us show that the third item implies the second one.  It follows
from the third item that either there are not directed triangles in $T$, or
every pair of directed triangles share an arc.   Suppose that $C_1, C_2$ and
$C_3$ are directed triangles in $T$ such that $C_1 = (u_1, u_2, u_3), C_2 =
(v_1, v_2, v_3)$, and $C_3 = (w_1, w_2, w_3)$.   Suppose without loss of
generality that $C_1$ and $C_2$ share the arc $(u_1, u_2) = (v_1, v_2)$ and
$C_2$ and $C_3$ share the arc $(v_2, v_3) = (w_2, w_3)$.   Now, $C_1$ and $C_3$
also share an arc.   It cannot be the case that $(u_1, u_2) = (w_2, w_3)$ or
$(u_1, u_2) = (w_3, w_1)$, as the former would imply that $v_1 = v_2$ and the
latter would imply that there is a loop at $u_1$. It is also impossible that
$(u_2, u_3) = (w_1, w_2)$ or $(u_2, u_3) = (w_3, w_1)$, as the former would
imply $u_2 = u_3$ and the latter would imply that $u_2 \to u_1$.   Finally, the
cases $(u_3, u_1) = (w_1, w_2)$ and $(u_3, u_1) = (w_2, w_3)$ are also absurd,
because we can deduce from both cases that there is a loop at $u_1$.   Hence, it
must be the case that $(u_1, u_2) = (w_1, w_2)$, $(u_2, u_3) = (w_2, w_3)$ or
$(u_3, u_1) = (w_3, w_1)$.   In the first two cases, we conclude that $C_1,
C_2$, and $C_3$ share an arc.   In the last case, $u_1 = v_1 = w_1$ and $u_3 =
v_3 = w_3$, so again, $C_1, C_2$, and $C_3$ share an arc.   We conclude that
every three directed triangles in $T$ share an arc, and thus, all directed
triangles in $T$ share an arc (use induction and the fact that if two directed
triangles share two arcs, then they are the same triangle).   Therefore, $T$ has
an arc such that its reversal results in a transitive tournament, proving the
second item.
\end{proof}

We now consider strongly chordal tournaments by allowing loops to be
present or absent. In a tournament $T$ with possible loops, we say a set of
vertices is {\em acyclic} if in $T$ it contains no directed cycle (other than a
loop).

The following lemma can be verified by a lengthy but straightforward
calculation.

\begin{lemma}
\label{lem:irr-gen}
Let $T$ be a tournament obtained from a tournament in the family $\cal{T}$ by
adding loops to an acyclic set of vertices, and such that the resulting
tournament does not contain $T_0$ (from Figure \ref{mixed_tournament}) as a
subgraph.   Then $T$ is a minimal obstruction for strong chordality.
\end{lemma}

Lemma \ref{lem:irr-gen} will be used multiple times in the proof of our
following theorem.   Let $\mathcal{T}^\circ$ be the family of tournaments
including $T_0$, the reflexive triangle $C_3$, and every possible tournament
obtained by a tournament in $\mathcal{T}$ by adding loops to a set of vertices
inducing a transitive tournament such that $T_0$ is not created.

\begin{theorem}
\label{thm:mix-char}
    Let $T$ be a tournament with possible loops.   The following statements are equivalent.
    \begin{enumerate}
        \item $T$ is strongly chordal.
        \item $T$ is obtained from a transitive tournament by reversing at most
            one arc with at least one irreflexive endvertex.
        \item $T$ is $\mathcal{T}^\circ$-free.
        \item $T$ is $T_0$-free, the irreflexive tournament obtained from $T$ by
            removing each loop is strongly chordal and the reflexive vertices of
            $T$ induce a transitive tournament.
    \end{enumerate}
\end{theorem}

\begin{proof}
Clearly, the third and fourth items are equivalent.   It also follows from
\Cref{lem:irr-gen} that the first item implies the third one.   To verify that
the fourth item implies the second one, notice that by removing all loops to
obtain an irreflexive tournament $T'$, \Cref{thm:irr-tou-cha} implies that $T'$
can be obtained by a transitive tournament $S$ by reversing at most one arc, let
us say $(u,v)$.   Suppose that both $u$ and $v$ are reflexive.   If there are
two vertices in the topological ordering of $S$ between $u$ and $v$, then either
any two of them are loopless, and together with $u$ and $v$ they induce a copy
of $T_0$, or one of them $w$ has a loop and $\{u, v, w\}$ induces a reflexive
$C_3$; both cases are impossible.   If there is only one vertex $w$ between $u$
and $v$ in $S$, by a similar argument as the one just used, we conclude that $w$
is an irreflexive vertex. Thus, $T'$ is also obtained from $S$ by reversing the
arc $(u,w)$, and such an arc has one irreflexive endpoint.

So, it remains to verify that the second item implies the first one.   We
consider two cases, when $T$ is strong and when it is non-strong.   When $T$ is
strong, notice that either the same $\Gamma$-free ordering used for $TT^\ast_n$
in the irreflexive case, or its reverse, will also work for this case.  The only
case, up to symmetry, where the order needs to be reversed, is when the vertices
are ordered by decreasing out-degree, and the first vertex has a loop.   Also,
if $n \ge 4$, then it can never happen that the first and last vertices are
reflexive, otherwise $T$ would contain $T_0$.

Notice that adding a source or a sink to a $\Gamma$-free tournament will result
again in a $\Gamma$-free tournament, regardless of whether the new vertex is
reflexive or irreflexive.   To obtain a $\Gamma$-free ordering for the new
tournament, it suffices to add the new vertex at the end of the previous
ordering. Thus, indeed the tournaments described in the theorem are strongly
chordal.

Let $T$ be a strong tournament which is strongly chordal.   If the underlying
irreflexive tournament $T^\circ$ of $T$ is isomorphic to $TT^\ast_n$, then $T$
does not contain $T_0$ as an induced subgraph, and it has the desired form.
Else, by Theorem \ref{thm:irr-tou-cha}, $T^\circ$ contains $T_i$ as a
subtournament, for $i \in \{ 1, \dots, 5 \}$.    If $T_i$ is also a tournament
of $T$, then $T$ is not strongly chordal, a contradiction.   Thus, $T$ contains
a copy of $T_i$ where some vertices are reflexive.   But this is not possible
either, because the directed reflexive triangle is a minimal obstruction for
strong chordality, as well as $T_0$ and $T_i$ with any acyclic subset of
vertices being reflexive, and not containing $T_0$. Thus, $T$ must have the
structure described in the first item of the theorem.

Now, if $T$ is non-strong, then every strong component of $T$ is either a single
vertex or contains a directed triangle.   Since the reflexive $3$-cycle and each
tournament obtained from $T_6$ by adding loops to an arbitrary acyclic subset
are minimal obstructions for strong chordality, it follows that at most one
connected component is not a single vertex.   Hence, the only non-trivial strong
component of $T$ has the structure described by the first item of this theorem,
and thus, $T$ has the desired structure.
\end{proof}

\begin{corollary}
    There is an $O(|V|+|A|)$ certifying algorithm to test whether a tournament
    with possible loops is strongly chordal.
\end{corollary}
\begin{proof}
    First, we calculate the strong components of the tournament in time
    $O(|V|+|A|)$.   If there is more than one non-trivial component, then the
    algorithm returns two disjoint directed triangles as a no-certificate. Now,
    consider the only non-trivial strong component $C$ (if any).   If $C$ has
    order $3$ or $4$ we can check whether it contains a reflexive $C_3$ or a
    copy of $T_0$, and if not, we can easily identify which is the reversed arc.
    Otherwise, let $k$ be the order of $C$.
    
    By looking at the degree sequence of $C$ we can identify which arc should be
    the reversed one: all out-degrees should be different, except for two pairs
    of vertices (the first and second having out degree $k-2$, and the last and
    second to last having out-degree $1$).   When this happens, a strong
    ordering can be obtained only using the out-degree sequence.   If this is
    not the case, choose a vertex of maximum out-degree $u$, suppose first that
    $u$ has in-degree at least $2$ and let $x$ and $y$ be two of its
    in-neighbours.  Since $T$ is strongly connected, and $u$ is a king in $T$
    there are vertices $x'$ and $y'$ such that $(u,x',x)$ and $(u,y',y)$ are
    paths in $T$.   If $x' \ne y'$, then $\{ u, x, x', y, y' \}$ induces a
    tournament with two arc-disjoint triangles, which is returned as a
    no-certificate.  If $x' = y'$, then, by the choice of $u$ there is a vertex
    $w$ which is an out-neighbour of $u$ but not of $x'$ (and hence it is an
    in-neighbour of $x'$).   If $(x,w)$ is an arc of $D$, then $(x',x,w,x')$ and
    $(x',y,u,x')$ are arc-disjoint triangles and otherwise $(u,w,x,u)$ and
    $(u,x',y,u)$ are arc-disjoint triangles; in either case $\{u, w, x, x', y\}$
    is returned as a no-certificate.   The algorithm proceeds analogously is the
    maximum in-degree of $T$ is at least $2$.

    Suppose now that the maximum out-degree and the maximum in-degree are both
    equal to $1$.   Choose a vertex $u$ of maximum out-degree and a vertex $v$
    of maximum in-degree.   If $(v,u)$ is an arc of $D$, then delete $u$ and $v$
    from $T$.   Since $T - \{ u, v \}$ is not transitive (because of the
    out-degree sequence), we can again find a strong component, and thus a
    triangle, in time $O(|V|+|A|)$.   Since $u$ dominates each vertex and $v$ is
    dominated by every vertex in such a triangle, it is clear that the vertices
    of the triangle together with $u$ and $v$ contain two arc-disjoint
    triangles, and thus can be returned as a no-certificate.   If the only
    in-neighbour of $u$ is equal to the only out-neighbour of $v$, call it $w$,
    again, $T-\{u,v\}$ is not transitive due to the degree sequence.  Hence, we
    can find a triangle in $T-\{u,v\}$ in time $O(|V|+|A|)$; if such a triangle
    uses $w$, then return the vertices of the triangle together with $u$ and $v$
    as a no-certificate.   Otherwise, return the vertices of the triangle
    together with $u, v$ and $w$ as a no-certificate.
    
    Finally, suppose that the only in-neigbour of $u$, call it $x$ is different
    to the only out-neighbour of $v$, call it $y$.   If $y$ dominates $x$, then
    choose any vertex $z$ different from $u, v, x$ and $y$.   When $y$ dominates
    $z$, then $(y,z,v,y)$ and $(y,x,u,y)$ are arc-disjoint triangles and $\{ u,
    v, x, y, z \}$ is returned as a no-certificate.   If for every possible
    choice of $z$, we have that $z$ dominates $y$, then $z$ also has out-degree
    equal to $1$, and we are in the case when the only in-neighbour of $u$
    coincides with the only out-neighbour of $v$.   When $z$ dominates $y$,
    since $T$ is strongly connected, $u$ is a king and $v$ is a serf, then there
    are vertices $w$ and $z$ such that $(u,w,x,w)$ and $(v,y,z,v)$ are
    arc-disjoint triangles, so $\{ u, v, w, x, y, z \}$ is returned as a
    no-certificate.

    Notice that finding each of the vertices described above takes only
    $O(|V|+|A|)$ time, and there are a constant number of searches for any
    output.
\end{proof}

\section{Irreflexive bipartite tournaments}
\label{sec:irr-bip}

Given a bipartite tournament $D$ with bipartition $(X,Y)$, its \textit{one-way
bigraphs} $B_X(D)$ and $B_Y(D)$ are the underlying graphs of the subdigraphs of
$D$ obtained from the deletion of every $YX$-arc, and every $XY$-arc,
respectively.   Clearly, for any bipartite tournament $D$, both $B_X(D)$ and
$B_Y(D)$ are bigraphs.   Let $M$ be the adjacency matrix of $D$ with some fixed
ordering of its vertex set, and call $M_X$ its submatrix induced by rows
corresponding to vertices in $X$ and columns corresponding to vertices in $Y$;
$M_Y$ is defined analogously.   Notice that $M_X$ and $M_Y$ correspond precisely
to the bi-adjacency matrices of $B_X(D)$ and $B_Y(D)$, respectively.   Moreover,
if we denote by $\overline{M}$ the matrix obtained from $M$ by swapping $0$'s
and $1$'s, it is clear that $M_X = \overline{M_Y^t}$.

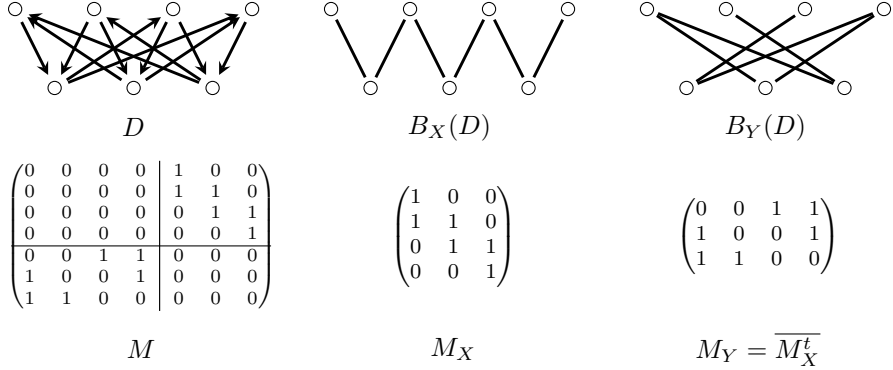
\begin{figure}[ht!]
    \begin{tikzpicture}
        \begin{scope}[scale=0.7]
            \node (1) [vertex] at (-2.25,0.75){};
            \node (2) [vertex] at (-0.75,0.75){};
            \node (3) [vertex] at (0.75,0.75){};
            \node (4) [vertex] at (2.25,0.75){};
            \node (5) [vertex] at (-1.5,-0.75){};
            \node (6) [vertex] at (0,-0.75){};
            \node (7) [vertex] at (1.5,-0.75){};
            \node (l) at (0,-1.5){$D$};

            \foreach \u/\v in {1/5,2/5,2/6,3/6,3/7,4/7}
                \draw [arc] (\u) to (\v);
            \foreach \u/\v in {5/3,5/4,6/1,6/4,7/1,7/2}
                \draw [arc] (\u) to (\v);
        \end{scope}
        \begin{scope}[scale=0.7,xshift=6cm]
            \node (1) [vertex] at (-2.25,0.75){};
            \node (2) [vertex] at (-0.75,0.75){};
            \node (3) [vertex] at (0.75,0.75){};
            \node (4) [vertex] at (2.25,0.75){};
            \node (5) [vertex] at (-1.5,-0.75){};
            \node (6) [vertex] at (0,-0.75){};
            \node (7) [vertex] at (1.5,-0.75){};
            \node (l) at (0,-1.5){$B_X(D)$};

            \foreach \u/\v in {1/5,2/5,2/6,3/6,3/7,4/7}
                \draw [edge] (\u) to (\v);
        \end{scope}
        \begin{scope}[scale=0.7,xshift=12cm]
            \node (1) [vertex] at (-2.25,0.75){};
            \node (2) [vertex] at (-0.75,0.75){};
            \node (3) [vertex] at (0.75,0.75){};
            \node (4) [vertex] at (2.25,0.75){};
            \node (5) [vertex] at (-1.5,-0.75){};
            \node (6) [vertex] at (0,-0.75){};
            \node (7) [vertex] at (1.5,-0.75){};
            \node (l) at (0,-1.5){$B_Y(D)$};

            \foreach \u/\v in {5/3,5/4,6/1,6/4,7/1,7/2}
                \draw [edge] (\u) to (\v);
        \end{scope}
        \begin{scope}[xshift=0.1cm,yshift=-2.5cm]
            \node (M) at (0,0){
            \scriptsize{
            $
            \begin{pNiceArray}{cccc|ccc}
                0 & 0 & 0 & 0 & 1 & 0 & 0 \\
                0 & 0 & 0 & 0 & 1 & 1 & 0 \\
                0 & 0 & 0 & 0 & 0 & 1 & 1 \\
                0 & 0 & 0 & 0 & 0 & 0 & 1 \\
                \hline
                0 & 0 & 1 & 1 & 0 & 0 & 0 \\
                1 & 0 & 0 & 1 & 0 & 0 & 0 \\
                1 & 1 & 0 & 0 & 0 & 0 & 0 \\
            \end{pNiceArray}
            $
            }
            };
            \node (l) at (0,-1.5){$M$};
        \end{scope}
        \begin{scope}[xshift=4.25cm,yshift=-2.5cm]
            \node (M) at (0,0){
            \footnotesize{
            $
            \begin{pNiceArray}{ccc}
                1 & 0 & 0 \\
                1 & 1 & 0 \\
                0 & 1 & 1 \\
                0 & 0 & 1 \\
            \end{pNiceArray}
            $
            }
            };
            \node (l) at (0,-1.5){$M_X$};
        \end{scope}
        \begin{scope}[xshift=8.3cm,yshift=-2.5cm]
            \node (M) at (0,0){
            \footnotesize{
            $
            \begin{pNiceArray}{cccc}
                0 & 0 & 1 & 1 \\
                1 & 0 & 0 & 1 \\
                1 & 1 & 0 & 0 \\
            \end{pNiceArray}
            $
            }
            };
            \node (l) at (0,-1.5){$M_Y = \overline{M_X^t}$};
        \end{scope}
    \end{tikzpicture}
    \caption{A bipartite tournament $D$, the one-way bigraphs $B_X(D)$ and
    $B_Y(D)$, and their corresponding adjacency or bi-adjacency matrices.}
    \label{fig:one-way-big}
\end{figure}

When $D$ is strongly chordal, the restrictions of the $\Gamma$-free ordering of
$D$ to the bigraphs $B_X(D)$ and $B_Y(D)$ are also (bipartite) $\Gamma$-free
orderings, so $B_X(D)$ and $B_Y(D)$ are strongly chordal bigraphs.   Conversely,
in order for a bipartite tournament $D$ to have a $\Gamma$-free ordering, it is
clear that we need $B_X(D)$ and $M_Y(D)$ to have $\Gamma$-free orderings when
considered as bigraphs, i.e., we want $M_X$ and $M_Y$ to have $\Gamma$-free
orderings, and while these orderings are not necessarily symmetric, the same
ordering of the vertices must work simultaneously for both matrices.   This
simultaneous $\Gamma$-free ordering condition prevents us from simply
considering bipartite tournaments where $B_X(D)$ and $B_Y(D)$ are both strongly
chordal.   In \Cref{fig:one-way-big} we see an example of a bipartite tournament
such that $B_X(D)$ and $B_Y(D)$ are both isomorphic to $P_7$, with the given
ordering $M_X$ is $\Gamma$-free, but $M_Y$ is not (actually, $D$ is not strongly
chordal).

Let $\mathcal{B}$ be the family of all bigraphs which are one-way bigraphs of a
strongly chordal bipartite tournament.   Observe that a characterization of
$\mathcal{B}$ directly yields a characterization of strongly chordal bipartite
tournaments.   In order to obtain such a characterization, we start by noticing
that $\mathcal{B}$ is a subclass of a well structured family of bigraphs. Since
we can interchange the roles of $X$ and $Y$ in a bipartition of a bipartite
tournament, we have that $\mathcal{B}$ is clearly closed under taking
\textit{bipartite complements} or \textit{bi-complements}, this is, the
complement of a bipartite graph relative to its complete bipartite spanning
supergraph.

\begin{figure}[ht!]
    \begin{tikzpicture}
        \begin{scope}[scale=0.8]
            \node (1) [vertex] at (-2.25,0.75){};
            \node (2) [vertex] at (-0.75,0.75){};
            \node (3) [vertex] at (0.75,0.75){};
            \node (4) [vertex] at (2.25,0.75){};
            \node (5) [vertex] at (-1.5,-0.75){};
            \node (6) [vertex] at (0,-0.75){};
            \node (7) [vertex] at (1.5,-0.75){};
            \node (l) at (0,-1.5){$B_1$};

            \foreach \u/\v in {1/5,2/5,2/6,3/6,3/7,4/7}
                \draw [arc] (\u) to (\v);
            \foreach \u/\v in {5/3,5/4,6/1,6/4,7/1,7/2}
                \draw [arc] (\u) to (\v);
        \end{scope}
        \begin{scope}[xshift=4.1cm,scale=0.8]
            \node (1) [vertex] at (-2.25,0.75){};
            \node (2) [vertex] at (-0.75,0.75){};
            \node (3) [vertex] at (0.75,0.75){};
            \node (4) [vertex] at (2.25,0.75){};
            \node (5) [vertex] at (-2.25,-0.75){};
            \node (6) [vertex] at (-0.75,-0.75){};
            \node (7) [vertex] at (0.75,-0.75){};
            \node (8) [vertex] at (2.25,-0.75){};
            \node (l) at (0,-1.5){$B_2$};

            \foreach \u/\v in {1/6,2/5,2/6,2/7,3/6,3/7,3/8,4/7}
                \draw [arc] (\u) to (\v);
            \foreach \u/\v in {5/1,5/3,5/4,6/4,7/1,8/1,8/2,8/4}
                \draw [arc] (\u) to (\v);
        \end{scope}
        \begin{scope}[xshift=8.2cm,scale=0.8]
            \node (1) [vertex] at (-2.25,0.75){};
            \node (2) [vertex] at (-0.75,0.75){};
            \node (3) [vertex] at (0.75,0.75){};
            \node (4) [vertex] at (2.25,0.75){};
            \node (5) [vertex] at (-1.5,-0.75){};
            \node (6) [vertex] at (0,-0.75){};
            \node (7) [vertex] at (1.5,-0.75){};
            \node (l) at (0,-1.5){$B_3$};

            \foreach \u/\v in {1/5,1/6,2/6,3/6,3/7,4/7}
                \draw [arc] (\u) to (\v);
            \foreach \u/\v in {5/2,5/3,5/4,6/4,7/1,7/2}
                \draw [arc] (\u) to (\v);
        \end{scope}
        \begin{scope}[yshift=-2.5cm,scale=0.8]
            \node (1) [vertex] at (-2.25,0.75){};
            \node (2) [vertex] at (-0.75,0.75){};
            \node (3) [vertex] at (0.75,0.75){};
            \node (4) [vertex] at (2.25,0.75){};
            \node (5) [vertex] at (-1.5,-0.75){};
            \node (6) [vertex] at (0,-0.75){};
            \node (7) [vertex] at (1.5,-0.75){};
            \node (l) at (0,-1.5){$P_7$};

            \foreach \u/\v in {1/5,2/5,2/6,3/6,3/7,4/7}
                \draw [edge] (\u) to (\v);
        \end{scope}
        \begin{scope}[xshift=4.1cm,yshift=-2.5cm,scale=0.8]
            \node (1) [vertex] at (-2.25,0.75){};
            \node (2) [vertex] at (-0.75,0.75){};
            \node (3) [vertex] at (0.75,0.75){};
            \node (4) [vertex] at (2.25,0.75){};
            \node (5) [vertex] at (-2.25,-0.75){};
            \node (6) [vertex] at (-0.75,-0.75){};
            \node (7) [vertex] at (0.75,-0.75){};
            \node (8) [vertex] at (2.25,-0.75){};
            \node (l) at (0,-1.5){$\textnormal{Sun}_4$};

            \foreach \u/\v in {1/6,2/5,2/6,2/7,3/6,3/7,3/8,4/7}
                \draw [edge] (\u) to (\v);
        \end{scope}
        \begin{scope}[xshift=8.2cm,yshift=-2.5cm,scale=0.8]
            \node (1) [vertex] at (-2.25,0.75){};
            \node (2) [vertex] at (-0.75,0.75){};
            \node (3) [vertex] at (0.75,0.75){};
            \node (4) [vertex] at (2.25,0.75){};
            \node (5) [vertex] at (-1.5,-0.75){};
            \node (6) [vertex] at (0,-0.75){};
            \node (7) [vertex] at (1.5,-0.75){};
            \node (l) at (0,-1.5){$S_{1,2,3}$};

            \foreach \u/\v in {1/5,1/6,2/6,3/6,3/7,4/7}
                \draw [edge] (\u) to (\v);
        \end{scope}
        \begin{scope}[yshift=-5cm,scale=0.8]
            \node (1) [vertex] at (-1.5,0.75){};
            \node (2) [vertex] at (0,0.75){};
            \node (3) [vertex] at (1.5,0.75){};
            \node (4) [vertex] at (-1.5,-0.75){};
            \node (5) [vertex] at (0,-0.75){};
            \node (6) [vertex] at (1.5,-0.75){};
            \node (l) at (0,-1.5){$C_6$};

            \foreach \u/\v in {4/2,4/3,5/1,5/3,6/1,6/2}
                \draw [edge] (\u) to (\v);
        \end{scope}
        \begin{scope}[xshift=4.1cm,yshift=-5cm,scale=0.8]
            \node (1) [vertex] at (-1.5,0.75){};
            \node (2) [vertex] at (0,0.75){};
            \node (3) [vertex] at (1.5,0.75){};
            \node (4) [vertex] at (-1.5,-0.75){};
            \node (5) [vertex] at (0,-0.75){};
            \node (6) [vertex] at (1.5,-0.75){};
            \node (l) at (0,-1.5){$B_4$};

            \foreach \u/\v in {1/4,2/5,3/6}
                \draw [arc] (\u) to (\v);
            \foreach \u/\v in {4/2,4/3,5/1,5/3,6/1,6/2}
                \draw [arc] (\u) to (\v);
        \end{scope}
        \begin{scope}[xshift=8.2cm,yshift=-5cm,scale=0.8]
            \node (1) [vertex] at (-1.5,0.75){};
            \node (2) [vertex] at (0,0.75){};
            \node (3) [vertex] at (1.5,0.75){};
            \node (4) [vertex] at (-1.5,-0.75){};
            \node (5) [vertex] at (0,-0.75){};
            \node (6) [vertex] at (1.5,-0.75){};
            \node (l) at (0,-1.5){$3K_2$};

            \foreach \u/\v in {1/4,2/5,3/6}
                \draw [edge] (\u) to (\v);
        \end{scope}
    \end{tikzpicture}
    \caption{Minimal obstructions for strong chordality in bipartite
    tournaments and their one-way bigraphs.}
    \label{fig:irr-bip-obs}
\end{figure}
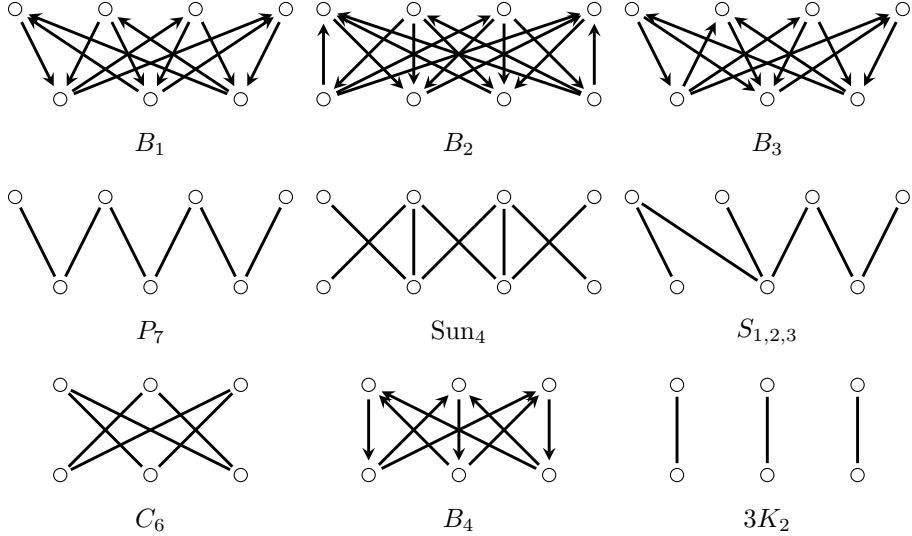

\begin{lemma}
\label{lem:irr-bip-obs}
    The bipartite tournaments $B_1, B_2, B_3$ and $B_4$, depicted in
    \Cref{fig:irr-bip-obs}, are minimal obstructions for strong chordality.
\end{lemma}
\begin{proof}
    It is easy to find to vertex disjoint copies of $C_4$ in $B_2$. Afterwards,
    the first item of \Cref{lem:irr-cyc} is enough to show that no vertex in
    $B_2$ is simple.   Similarly, it is easy to find three copies of $C_4$
    pairwise sharing an arc in $B_4$; for each of such cycles, the two vertices
    not in the cycle have two in-neighbours or two out-neighbours in the cycle,
    so again we deduce from \Cref{lem:irr-cyc} that no vertex in $B_4$ is
    simple.

    Suppose that $(X,Y)$ is a bipartition of $B_1$ with $|X| = 4$.   By using
    the first item of \Cref{lem:irr-cyc} with suitable choices of $4$-cycles, it
    is possible to conclude that neither vertices in $X$ nor the vertex in $Y$
    having both vertices of out-degree $1$ as out-neighbours is simple.   For
    the remaining vertices in $Y$, call them $u$ and $v$, notice that the
    in-neighbourhood of one is the out-neighbourhood of the other and vice
    versa.   Moreover, \Cref{lem:simple-ord} imposes contradictory orderings on
    their in- and out-neighbourhoods, so starting an ordering with $u$ implies
    that at least one of its in-neighbours and one of its out-neighbours must be
    chosen as the next vertex in the ordering before chosing $v$.   Nonetheless,
    the second vertex in a possible $\Gamma$-free ordering of $B_1$ must be the
    third vertex in $Y$, as the vertices in $X$ that are permitted based on the
    restrictions of \Cref{lem:simple-ord} are not simple.   But again, choosing
    such a vertex imposes restrictions through \Cref{lem:simple-ord} on vertices
    in $X$ which are incompatible with the previous ones.   Therefore, after
    choosing two vertices, it is impossible to complete a $\Gamma$-free ordering
    for $B_1$.   An analogous argument shows that $B_3$ is not strongly chordal.
\end{proof}

Giakoumakis and Vanherpe introduced \textit{co-bigraphs} in
\cite{giakoumakisAAM18}, a bipartite analogue for cographs defined as the
bigraphs that can be reduced to isolated vertices by recursively taking
bi-complements of their connected components.   They proved that a bigraph is a
co-bigraph if and only if it is $\{P_7, \textnormal{Sun}_4, S_{1,2,3}\}$-free.
It follows from \Cref{lem:irr-bip-obs} and the characterization of co-bigraphs
that both one-way bigraphs of an irreflexive bipartite tournament are
co-bigraphs. We use this information, together with the fact that $B_4$ is a
minimal obstruction for strong chordality, to obtain the following necessary
condition for an irreflexive bipartite tournament to be strongly chordal. Recall
that a graph is a \textit{bipartite chain graph} if it is bipartite and the
neighbourhoods of the vertices in each part can be linearly ordered with respect
to inclusion. It is well-known that the complete list of minimal obstructions
for the family of bipartite chain graphs is $\{C_3, 2K_2, C_5\}$.

Our next results takes advantage of the structure of co-bigraphs and bipartite
chain graphs.

\begin{lemma}
\label{lem:irr-chain}
    Let $D$ be an irreflexive bipartite tournament.   If $D$ is strongly
    chordal, then one of $B_X(D)$ or $B_Y(D)$ is a disconnected graph with at
    most two non-trivial components, each of which is a bipartite chain graph.
\end{lemma}
\begin{proof}
    As we observed in the above paragraph, strongly one-way bigraphs of chordal
    bipartite tournaments are co-bigraphs.    Since co-bigraphs are precisely
    the bipartite graphs which are reducible to isolated vertices by applying
    the bipartite complement operation, and the bipartite complement of $B_X(D)$
    is $B_Y(D)$ (and vice versa), it follows that one of $B_X(D)$ and $B_Y(D)$
    is disconnected.   Suppose, without loss of generality, that $B_X(D)$ is
    disconnected.   Given that $B_4$ is a minimal obstruction for strong
    chordality, we have that $B_X(D)$ is $3K_2$-free, and hence, it has at most
    two non-trivial connected components.   But again, the fact that $B_X(D)$ is
    $3K_2$-free implies that such components are $2K_2$-free, and hence they are
    bipartite chain graphs.
\end{proof}

The main result of this section contains characterizations of a strongly chordal
irreflexive bipartite tournament $D$ in terms of forbidden induced subdigraphs,
the structure of the graphs $B_X(D)$ and $B_Y(D)$, and the global structure of
$D$.

\begin{theorem}
\label{thm:irr-bip-char}
    Let $D$ be an irreflexive bipartite tournament.   The following statements
    are equivalent.
    \begin{enumerate}
        \item $D$ is strongly chordal.
        \item $D$ is $\{B_1, B_2, B_3, B_4\}$-free.
        \item $B_X(D)$ and $B_Y(D)$ are $\{3K_2, C_6\}$-free co-bigraphs.
        \item $B_X(D)$ or $B_Y(D)$ is a disconnected graph with at most two
            non-trivial connected components, each of which is a bipartite chain
            graph.
    \end{enumerate}
\end{theorem}
\begin{proof}
    \Cref{lem:irr-bip-obs} proves that the first item implies the second.   The
    discussion before \Cref{lem:irr-chain} and the proof of the same lemma show
    that the second item implies the third and that the third implies the
    fourth. To conclude the proof, it suffices to prove that the last item
    implies the first one.

    Suppose that $B_X(D)$ is a disconnected graph with non-trivial components
    $C_1$ and $C_2$ having bipartitions $(X_1, Y_1)$ and $(X_2, Y_2)$,
    respectively.   Suppose also that $x_{i1}, \dots, x_{ir_i}$ and $y_{i1},
    \dots, y_{is_i}$ are linear orderings of $X_i$ and $Y_i$, respectively,
    satisfying $N(x_j) \subseteq N(x_k)$ if $j < k$ and $N(y_j) \subseteq
    N(y_k)$ if $j < k$, for $i \in \{1,2\}$.   It is a routine exercise to
    verify that the ordering obtained by starting with the vertices of $X_1$
    followed by $Y_2, Y_1$ and $X_2$, such that the first two are in the given
    linear oderings and the last two follow the reverse of such orderings, is a
    $\Gamma$-free ordering of $D$.
\end{proof}

We finalize this section by using \Cref{thm:irr-bip-char} to show that strongly
chordal irreflexive bipartite tournaments can be recognized in linear time. Our
algorithm is not certifying, but it fails to return a no-certificate in only one
output.   We use the bipartite chain graph certifying recognition algorithm of
Heggernes and Kratsch \cite{heggernessNJC14} as a subroutine.

\begin{corollary}
\label{cor:irr-bip-alg}
    There is an $O(|V|+|A|)$-time algorithm to recognize strongly chordal
    irreflexive bipartite tournaments.
\end{corollary}
\begin{proof}
    For an input irreflexive bipartite tournament $D$ with bipartition $(X,Y)$,
    we can build $B_X(D)$ and $B_Y(D)$ and find their connected components in
    time $O(|V|+|A|)$.  If both of them are connected, then they are not
    co-bigraphs, and the algorithm answers that $D$ is not strongly chordal
    (this is the only case where a certificate is not provided).   Else, we
    assume without loss of generality that $B_X(D)$ is disconnected.   If there
    are more than three non-trivial components, the algorithm returns $3K_2$ as
    a no-certificate. We use again $O(|V|+|A|)$ time to test each of the (at
    most) two connected components of $B_X(D)$ to be bipartite chain graphs.
    If they are not, a $2K_2$ is found in one of the components, and the
    algorithm returns $3K_2$ as a no-certificate.   At this point we know that
    $B_X(D)$ has exactly two connected components, and both are bipartite chain
    graphs.   Finally, the algorithm returns the ordering constructed in
    \Cref{thm:irr-bip-char} as a yes-certificate.
\end{proof}

\section{Irreflexive tournaments minus one arc}
\label{sec:irr-minus}

\begin{figure}[htb!]
    \begin{tikzpicture}
    \begin{scope}[xshift=0cm,scale=0.45,rotate=30]
        \foreach \i in {0,...,5}
            \draw ({(360/6)*\i}:2) node[vertex](\i){};
        \foreach \i/\j in {0/2,0/3,0/5,1/3,2/1,3/2,4/0,4/1,4/2,4/3,5/1,5/2,5/3,5/4}
            \draw [arc] (\i) to (\j);
        \node (l) at (95:2.5){$S_{15}$};
    \end{scope}
    \begin{scope}[xshift=2.5cm,scale=0.45,rotate=30]
        \foreach \i in {0,...,5}
            \draw ({(360/6)*\i}:2) node[vertex](\i){};
        \foreach \i/\j in {0/2,0/3,0/4,0/5,1/2,1/3,2/4,3/2,4/1,4/3,5/1,5/2,5/3,5/4}
            \draw [arc] (\i) to (\j);
        \node (l) at (95:2.5){$S_{16}$};
    \end{scope}
    \begin{scope}[xshift=5cm,scale=0.45,rotate=30]
        \foreach \i in {0,...,5}
            \draw ({(360/6)*\i}:2) node[vertex](\i){};
        \foreach \i/\j in {0/2,0/3,0/4,0/5,1/2,1/3,1/4,2/4,3/2,4/3,5/1,5/2,5/3,5/4}
            \draw [arc] (\i) to (\j);
        \node (l) at (95:2.5){$S_{17}$};
    \end{scope}
    \begin{scope}[xshift=7.5cm,scale=0.45,rotate=30]
        \foreach \i in {0,...,5}
            \draw ({(360/6)*\i}:2) node[vertex](\i){};
        \foreach \i/\j in {0/3,2/0,2/1,2/3,2/5,3/1,4/0,4/1,4/2,4/3,5/0,5/1,5/3,5/4}
            \draw [arc] (\i) to (\j);
        \node (l) at (95:2.5){$S_{18}$};
    \end{scope}
    \begin{scope}[xshift=10cm,scale=0.45,rotate=30]
        \foreach \i in {0,...,5}
            \draw ({(360/6)*\i}:2) node[vertex](\i){};
        \foreach \i/\j in {0/2,0/3,2/1,2/3,2/5,3/1,4/0,4/1,4/2,4/3,5/0,5/1,5/3,5/4}
            \draw [arc] (\i) to (\j);
        \node (l) at (95:2.5){$S_{19}$};
    \end{scope}
    \end{tikzpicture}
    \caption{Minimal obstructions for strongly chordal tournaments minus one arc on
    six vertices.}
    \label{fig:irr-minus-6}
\end{figure}
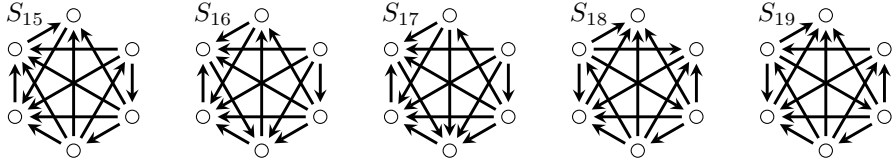

After characterizing bipartite tournaments it seems a natural next step to
consider multipartite tournaments.   Nonetheless, as the restricted case which
we develop in this section shows, the general case appears to be a lot more
complex than the tournament and bipartite tournament cases, at least regarding
the number and structure of minimal obstructions.   Together with the
tournaments shown in \Cref{fig:minobsirrtour} the multipartite tournaments in
\Cref{fig:irr-minus-5,fig:irr-minus-6} complete the minimal obstructions to
strong chordality for the family of tournaments minus an arc (or equivalently,
complete multipartite tournaments where each part has one vertex, except maybe
for one part that may have two vertices).   Let $\mathcal{S}$ be the family
$\mathcal{S} = \{ S_1, \dots, S_{19} \}$.

\begin{figure}[htb!]
\begin{tikzpicture}
    \begin{scope}[xshift=0cm,scale=0.5,rotate=18]
        \foreach \i in {0,...,4}
            \draw ({(360/5)*\i}:2) node[vertex](\i){};
        \foreach \i/\j in {0/4,1/3,2/0,2/1,3/0,3/2,4/1,4/2,4/3}
            \draw [arc] (\i) to (\j);
        \node (l) at (108:2.3){$S_1$};
    \end{scope}
    \begin{scope}[xshift=2.5cm,scale=0.5,rotate=18]
        \foreach \i in {0,...,4}
            \draw ({(360/5)*\i}:2) node[vertex](\i){};
        \foreach \i/\j in {0/2,0/4,1/3,2/1,3/0,3/2,4/1,4/2,4/3}
            \draw [arc] (\i) to (\j);
        \node (l) at (108:2.3){$S_2$};
    \end{scope}
    \begin{scope}[xshift=5cm,scale=0.5,rotate=18]
        \foreach \i in {0,...,4}
            \draw ({(360/5)*\i}:2) node[vertex](\i){};
        \foreach \i/\j in {0/3,0/4,1/3,2/0,2/1,3/2,4/1,4/2,4/3}
            \draw [arc] (\i) to (\j);
        \node (l) at (108:2.3){$S_3$};
    \end{scope}
    \begin{scope}[xshift=7.5cm,scale=0.5,rotate=18]
        \foreach \i in {0,...,4}
            \draw ({(360/5)*\i}:2) node[vertex](\i){};
        \foreach \i/\j in {0/2,0/4,1/4,2/1,3/0,3/1,3/2,4/2,4/3}
            \draw [arc] (\i) to (\j);
        \node (l) at (108:2.3){$S_4$};
    \end{scope}
    \begin{scope}[xshift=10cm,scale=0.5,rotate=18]
        \foreach \i in {0,...,4}
            \draw ({(360/5)*\i}:2) node[vertex](\i){};
        \foreach \i/\j in {0/3,0/4,1/4,2/0,2/1,3/1,3/2,4/2,4/3}
            \draw [arc] (\i) to (\j);
        \node (l) at (108:2.3){$S_5$};
    \end{scope}
    \begin{scope}[xshift=0cm,yshift=-2.5cm,scale=0.5,rotate=18]
        \foreach \i in {0,...,4}
            \draw ({(360/5)*\i}:2) node[vertex](\i){};
        \foreach \i/\j in {0/3,0/4,1/2,1/4,2/0,3/1,3/2,4/2,4/3}
            \draw [arc] (\i) to (\j);
        \node (l) at (108:2.3){$S_6$};
    \end{scope}
    \begin{scope}[xshift=2.5cm,yshift=-2.5cm,scale=0.5,rotate=18]
        \foreach \i in {0,...,4}
            \draw ({(360/5)*\i}:2) node[vertex](\i){};
        \foreach \i/\j in {0/2,0/3,2/1,2/4,3/1,3/2,4/0,4/1,4/3}
            \draw [arc] (\i) to (\j);
        \node (l) at (108:2.3){$S_7$};
    \end{scope}
    \begin{scope}[xshift=5cm,yshift=-2.5cm,scale=0.5,rotate=18]
        \foreach \i in {0,...,4}
            \draw ({(360/5)*\i}:2) node[vertex](\i){};
        \foreach \i/\j in {0/2,0/3,0/4,2/1,2/4,3/1,3/2,4/1,4/3}
            \draw [arc] (\i) to (\j);
        \node (l) at (108:2.3){$S_8$};
    \end{scope}
    \begin{scope}[xshift=7.5cm,yshift=-2.5cm,scale=0.5,rotate=18]
        \foreach \i in {0,...,4}
            \draw ({(360/5)*\i}:2) node[vertex](\i){};
        \foreach \i/\j in {0/3,1/2,2/0,2/4,3/1,3/2,4/0,4/1,4/3}
            \draw [arc] (\i) to (\j);
        \node (l) at (108:2.3){$S_9$};
    \end{scope}
    \begin{scope}[xshift=10cm,yshift=-2.5cm,scale=0.5,rotate=18]
        \foreach \i in {0,...,4}
            \draw ({(360/5)*\i}:2) node[vertex](\i){};
        \foreach \i/\j in {0/2,0/3,1/2,2/4,3/1,3/2,4/0,4/1,4/3}
            \draw [arc] (\i) to (\j);
        \node (l) at (108:2.3){$S_{10}$};
    \end{scope}
    \begin{scope}[xshift=1.25cm,yshift=-5cm,scale=0.5,rotate=18]
        \foreach \i in {0,...,4}
            \draw ({(360/5)*\i}:2) node[vertex](\i){};
        \foreach \i/\j in {0/2,0/4,1/2,2/4,3/0,3/1,3/2,4/1,4/3}
            \draw [arc] (\i) to (\j);
        \node (l) at (108:2.3){$S_{11}$};
    \end{scope}
    \begin{scope}[xshift=3.75cm,yshift=-5cm,scale=0.5,rotate=18]
        \foreach \i in {0,...,4}
            \draw ({(360/5)*\i}:2) node[vertex](\i){};
        \foreach \i/\j in {0/3,0/4,1/2,2/0,2/4,3/1,3/2,4/1,4/3}
            \draw [arc] (\i) to (\j);
        \node (l) at (108:2.3){$S_{12}$};
    \end{scope}
    \begin{scope}[xshift=6.25cm,yshift=-5cm,scale=0.5,rotate=18]
        \foreach \i in {0,...,4}
            \draw ({(360/5)*\i}:2) node[vertex](\i){};
        \foreach \i/\j in {0/2,0/3,0/4,1/2,2/4,3/1,3/2,4/1,4/3}
            \draw [arc] (\i) to (\j);
        \node (l) at (108:2.3){$S_{13}$};
    \end{scope}
    \begin{scope}[xshift=8.75cm,yshift=-5cm,scale=0.5,rotate=18]
        \foreach \i in {0,...,4}
            \draw ({(360/5)*\i}:2) node[vertex](\i){};
        \foreach \i/\j in {0/2,0/4,1/2,1/3,2/4,3/0,3/2,4/1,4/3}
            \draw [arc] (\i) to (\j);
        \node (l) at (108:2.3){$S_{14}$};
    \end{scope}
\end{tikzpicture}
\caption{Minimal obstructions for strongly chordal tournaments minus one arc on
five vertices.}
\label{fig:irr-minus-5}
\end{figure}
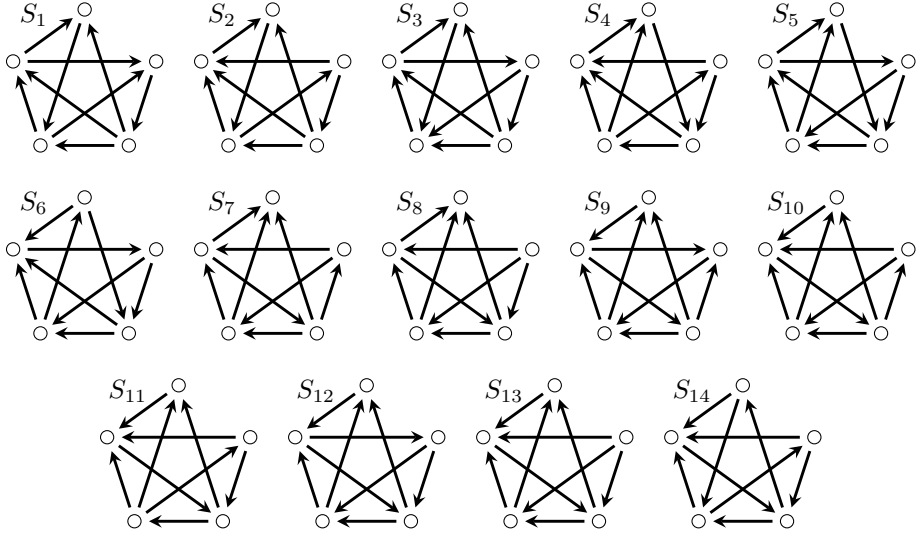

It is useful to explore the possible $\Gamma$-free orderings of an irreflexive
strongly chordal tournament.   If $T$ is transitive with topological ordering
$v_1 \to \cdots \to v_n$, then any ordering starting at $v_i$ and respecting the
topological ordering to the left and to the right of $v_i$, possibly alternating
between in-neighbours and out-neighbours of $v_i$, is a $\Gamma$-free ordering.
Moreover, by virtue of \Cref{lem:simple-ord} every $\Gamma$-free ordering of $T$
is of this form.   If $T$ is obtained from a transitive tournament with
topological ordering $v_1 \to \cdots \to v_n$ by reversing the arc $(v_i, v_j)$
with $i+1 < j$, then it follows from \Cref{lem:irr-cyc} that no vertex before
$v_i$ or after $v_j$ in the topological ordering is simple.   If $j = i+2$, then
any ordering placing the vertices in the only non-trivial strong component of
$T$ first, followed by all vertices to its left and then all vertices on its
right, respecting the topological ordering from it, is a $\Gamma$-free ordering.
Else, $i+2 < j$ and for any $k$ and $l$ such that $i < k < l < j$, we have that
$v_i$ and $v_k$ are in-neighbours of $v_l$ such that $v_k \in N^+(v_i) \setminus
N^+(v_k)$ and $v_j \in N^+(v_k) \setminus N^+(v_i)$, i.e., $v_i$ and $v_k$ have
incomparable out-neighbourhoods.   Hence, $v_l$ is not simple, and analogously,
$v_k$ is not simple either.   So, the only simple vertices are $v_i$ and $v_j$,
so a $\Gamma$-free ordering must start in one of those vertices.   Once one of
these two vertices is chosen to go first, the rest of the ordering is determined
by \Cref{lem:simple-ord} (which basically tell us to respect the topological
ordering of the rest of $T$).

\begin{lemma}
\label{lem:irr-minus-min-obs}
    All irreflexive multipartite tournaments in the family $\mathcal{S}$ are
    minimal obstructions to strong chordality.
\end{lemma}
\begin{proof}
    Non-adjacent vertices in $S_7, S_8, S_{10}$ and $S_{13}$ have two
    in-neighbours or two out-neighbours in a directed triangle, so the first
    item of \Cref{lem:irr-cyc} implies that they are not simple.   The second
    item of the same lemma implies that none of the remaining three vertices
    (forming a directed triangle) are simple.

    Using \Cref{lem:irr-cyc} it is easy to observe that there is no simple
    vertex in $S_1, S_5, S_6, S_9, S_{15}, S_{16}, S_{17}$, $S_{18}$ $S_{19}$,
    and also that there is only one simple non-peak vertex in $S_3, S_{11},
    S_{12}, S_{14}$.   A simple exploration shows that all vertices in $S_2$ are
    peak vertices. To see that $S_4$ is not strongly chordal, label its five
    vertices by $a, b, c, d, e$ in the clockwise order starting at the top
    vertex. We observe that $a, c$ are are the only simple vertices and $a$ is
    the only non-peak in $S_4$. So any strong ordering $<$ of $S_4$ must begin
    with $c$ and end with $a$. But then we have $b < a$ and $c < e$ for which
    $(b,c), (b,e), (a,c) \in E(S_4)$ and $(a,e) \notin E(S_4)$. This contradicts
    that $<$ is a strong ordering of $S_4$.    A similar argument proves that
    $S_3$ is not strongly chordal, but now (with the same labelling) $b$ and $d$
    are the only simple vertices and $d$ is the only non-peak vertex.

    The minimality of $S_i$ for $i \in \{ 1, \dots, 14 \}$ is easy to verify as
    every tournament on four vertices and each orientation of the diamond are
    strongly chordal.   Finally, $\Gamma$-free orderings for each of the
    remaining tournamentns when any vertex is deleted are not hard to obtain by
    following the description of the $\Gamma$-free orderings for strongly
    chordal tournaments in the paragraph previous to this lemma.
\end{proof}

Let $T$ be a tournament obtained from a transitive tournament with topological
ordering $1 \to \dots \to n$ by reversing the arc $(a,b)$ (so $a < b$ and
$(b,a)$ is an arc of $T$).   An arc $(x,y)$ of $T$ is \textit{critical} if any
of the following conditions hold:
\begin{enumerate}
    \item $x+1 < y < a$
    \item $x < a-1 < y-1 < b-2$
    \item $x < a < y-1 < b-1$
    \item $x < a < b-1 < y-1$
    \item $a < x < y-1 < b-1$
    \item $a < x < b-1 < y-1$
    \item $a+1 < x < b < y-1$
    \item $b < x < y-1$
\end{enumerate}

\begin{lemma}
\label{lem:irr-critical}
    Let $D$ be obtained from an irreflexive strongly chordal tournament by
    deleting an arc.   If the deleted arc is critical, then $D$ is not strongly
    chordal.
\end{lemma}
\begin{proof}
    We consider the deletion of each possible type of critical arc:
    \begin{enumerate}
        \item The digraph $S_{17}$ is an induced subdigraph of $D$.
        \item The digraph $S_{16}$ is an induced subdigraph of $D$.
        \item The digraph $S_{13}$ is an induced subdigraph of $D$.
        \item The digraph $S_8$ is an induced subdigraph of $D$.
        \item The digraph $S_{10}$ is an induced subdigraph of $D$.
        \item The digraph $S_7$ is an induced subdigraph of $D$.
        \item The digraph $S_{19}$ is an induced subdigraph of $D$.
        \item The digraph $S_{18}$ is an induced subdigraph of $D$.
    \end{enumerate}
\end{proof}

Let $T$ be a strongly chordal irreflexive tournament. By Theorem
\ref{thm:irr-tou-cha} $T$ is either a transitive tournament or is obtained from
a transitive tournament by reversing an arc. We call an arc in $T$ {\em
non-critical} if it is not critical. Note that, if $T$ is transitive, then every
arc in $T$ is non-critical and, if $T$ is obtained from a transitive tournament
by reversing an arc, then the reversed arc is non-critical. 

\begin{theorem}
\label{thm:irr-minus-char}
   Let $D$ be obtained from an irreflexive tournament $T$ by deleting the arc
   $(x,y)$. The following statements are equivalent.
    \begin{enumerate}
        \item $D$ is strongly chordal.
        \item $D$ is $(\mathcal{T} \cup \mathcal{S})$-free.
        \item $T$ strongly chordal and $(x,y)$ is non-critical or $T-y$ is
              strongly chordal and $y$ is a false twin of $x$ in $D$.
    \end{enumerate}
\end{theorem}
\begin{proof}
It follows from \Cref{pro:irr-min-obs,lem:irr-minus-min-obs} that the first item
implies the second one.   To prove that the third item implies the first one,
suppose that $T$ is strongly chordal. Then it has a strong ordering. One can
verify that the same ordering is a strong ordering of $D$. On the other hand
suppose that $T-y$ is strongly chordal and $y$ is a false twin of $x$ in $D$.
Then the vertex ordering of $D$ obtained from a strong ordering of $T - y$ by
inserting $y$ next to $x$ is a strong ordering of $D$.

It remains to prove that the second item implies the third one. So assume that
$D$ is $(\mathcal{T} \cup \mathcal{S})$-free. If $D$ is acyclic, then $T = D +
(x,y)$ is strongly chordal and $(x,y)$ is non-critical. Thus Statement~3 holds.
Since $D$ is $\mathcal{T}$-free, $T-y$ ($=D-y$) is also $\mathcal{T}$-free and
hence is strongly chordal by \Cref{thm:irr-tou-cha}. If $x$ and $y$ are false
twin in $D$ then Statement~3 holds. So from now on we assume that $D$ contains
directed cycles and that $x, y$ are not not false twins.

Since $D$ contains directed cycles, it contains directed triangles. We prove by
contradiction that $D$ does not contain two arc-disjoint directed triangles. So
suppose that $\Delta$ and $\Delta'$ are two arc-disjoint directed triangles in
$D$. Denote by $H$ the subdigraph of $D$ induced by $\Delta \cup \Delta'$. We
claim that $H$ is not a tournament (i.e., it contains both $x, y$). Indeed,
since $D$ is $\mathcal{T}$-free, $H$ is also $\mathcal{T}$-free. If $H$ is a
tournament then by \Cref{thm:irr-tou-cha} it does not contain two arc-disjoint
directed triangles, which contradicts the fact that $\Delta$ and $\Delta'$ are
arc-disjoint directed triangles in $H$.

The subdigraph $H$ has either five or six vertices depending on whether or not
$\Delta$ and $\Delta'$ have a common vertex. Suppose that $H$ has five vertices,
that is, $\Delta$ and $\Delta'$ have a common vertex. Denote $\Delta = abc$ and
$\Delta' = cfg$ where $c$ is the common vertex between the two triangles. As $H$
is not a tournament, it has exactly three arcs between $\{a,b\}$ and $\{f,g\}$.
By relabeling the vertices if necessary we can assume that $H$ has at most one
arc from $\{f,g\}$ to $\{a,b\}$. If $H$ has no arc from $\{f,g\}$ to $\{a,b\}$,
then it is $S_2$ when $a$ and $g$ are not adjacent, $S_3$ when $a$ and $f$ are
not adjacent, $S_4$ when $b$ and $g$ are not adjacent, and $S_{11}$ when $b$ and
$f$ are not adjacent. This contradicts the assumption that $D$ is
$\mathcal{S}$-free.

Suppose that $(g,a)$ is an arc. Then $H$ is $S_1$ when $a$ and $f$ are not
adjacent, $S_5$ when $b$ and $f$ are not adjacent, and $S_6$ when $b$ and $g$
are not adjacent. This again contradicts the assumption that $D$ is
$\mathcal{S}$-free. Suppose that $(g,b)$ is an arc. Then $H$ is $S_6$ when $a$
and $g$ are not adjacent and is $S_9$ when $f$ is not adjacent to $a$ or $b$.
Suppose that $(f,a)$ is an arc. Then $H$ is $S_1$ when $a$ and $f$ are not
adjacent, $S_5$ when $b$ and $g$ are not adjacent, and $S_{14}$ when $b$ and $f$
are not adjacent. Finally, suppose that $(f,b)$ is an arc. Then $H$ is $S_{12}$
when $a$ and $g$ are not adjacent. When $a$ and $f$ are not adjacent, they are
false twins in $H$. By our assumption above, $a, f$ (which are $x, y$) are not
false twins in $D$. So there is a vertex $w$ of which one of $a$ and $f$ is an
in-neighbour and the other is an out-neighbour. Assume without loss of
generality that $f$ is an in-neighbour and $a$ is an out-neighbour of of $w$. We
claim that $c, g$ are out-neighbours of $w$. Indeed if $g$ is an in-neighbour of
$w$, then $H-f+w$ is a tournament containing the arc-disjoint directed triangles
$wag$ and $abc$, which implies by \Cref{thm:irr-tou-cha} that $D$ is not
$\mathcal{T}$-free, a contradiction to the assumption. So $g$ is an
out-neighbour of $w$. If $c$ is an in-neighbour of $w$, then $H-a+w$ is a
tournament containing arc-disjoint directed triangles $wgc$ and $abc$, resulting
in a similar contradiction. So $c$ is an out-neighbour of $w$. We see now that
$H-g+w$ contains arc-disjoint directed triangles $wcf$ and $abc$ and $a, f$ are
not false twins in it. The above proof applied to $H-g+w$ (instead of $H$) leads
to a contradiction. The discussion for the case when $b$ and $g$ are not
adjacent is similar (in this case $b, g$ are false twins in $H$). 

Suppose now that $H$ has six vertices, that is, $\Delta$ and $\Delta'$ have no
common vertex. From the above proof we can assume that no two arc-disjoint
directed triangles in $D$ have a common vertex. It follows that the arcs between
$\Delta$ and $\Delta'$ are either from $\Delta$ to $\Delta'$ or from $\Delta'$
to $\Delta$. In either case $H$ is $S_{15}$, a contradiction to the assumption
that $D$ is $\mathcal{S}$-free. Hence $D$ does not contain arc-disjoint cycles.
This means that there exists an arc $(s,t)$ such that $D - (s,t)$ is acyclic. 

Let $1, \dots, n$ be a topological ordering of $D - (s,t)$. Then $s > t$. We may
assume that $x < y$. Clearly, $T = D + (x,y)$ is strongly chordal. Since $D$ is
$\mathcal{S}$-free, by \Cref{lem:irr-critical} $(x,y)$ must be a non-critical
arc in $T$. Therefore $D$ is obtained from the strongly chordal tournament $T$
by deleting the non-critical arc $(x,y)$. This completes the proof. 
\end{proof}

\begin{corollary}
\label{lem:irr-minus-two-cyc}
    Let $D$ be obtained from an irreflexive tournament by deleting the arc
    $(x,y)$.   If $x$ and $y$ are not twins in $D$ and $D$ contains two
    arc-disjoint triangles, then $D$ is not strongly chordal.
\end{corollary}
\begin{proof}
This follows from the proof of Theorem \ref{thm:irr-minus-char}.
\end{proof}

\section{Reflexive multipartite tournaments}
\label{sec:ref-mult}

Like reflexive strongly chordal tournaments, reflexive strongly chordal
multipartite tournaments also have a restricted structure.   Our first result of
the section deals with the underlying graph of reflexive strongly chordal
multipartite tournaments.

\begin{lemma}
\label{lem:ref-comsplit}
    Let $D$ be a reflexive multipartite tournament.   If $D$ is strongly
    chordal, then its underlying graph is a complete split graph.
\end{lemma}

\begin{proof}
    Being $D$ reflexive, the fact that it is strongly chordal implies that its
    underlying graph $U_D$ is chordal, so it is $C_4$-free.   But $U_D$ is also
    a complete multipartite graph, so it also $(K_1 + K_2)$-free.   Therefore,
    it is $\{ K_1 + K_2, C_4 \}$-free, i.e., a complete split graph.
\end{proof}

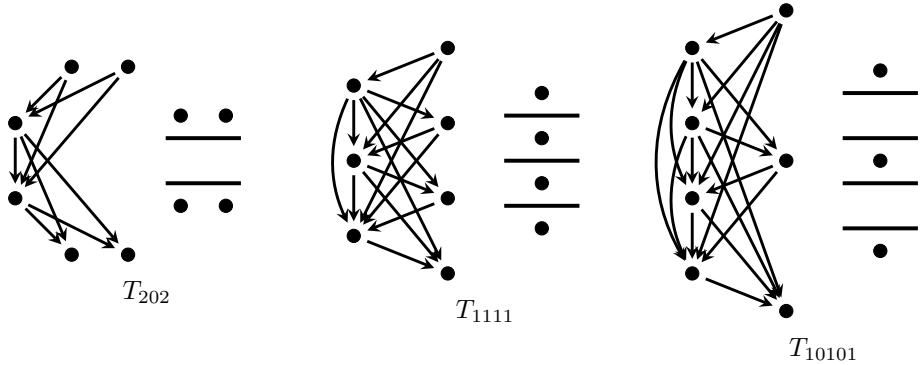
\begin{figure}[ht!]
\centering
\begin{tikzpicture}
    \begin{scope}
        \node (t1) [blackV] at (0,0.5){};
        \node (t2) [blackV] at (0,-0.5){};
        \node (s1) [blackV] at (0.75,1.25){};
        \node (s2) [blackV] at (1.5,1.25){};
        \node (s3) [blackV] at (0.75,-1.25){};
        \node (s4) [blackV] at (1.5,-1.25){};
        \node (l1) at (1.75,-1.75){$T_{202}$};

        \draw [arc] (t1) to (t2);
        \foreach \i in {1,2}{
            \draw [arc] (s1) to (t\i);
            \draw [arc] (s2) to (t\i);
            \draw [arc] (t\i) to (s3);
            \draw [arc] (t\i) to (s4);
        }
    \end{scope}
    \begin{scope}[xshift=2.5cm]
        \node (n1) [blackV] at (0.3,0.6){};
        \node (n1) [blackV] at (-0.3,0.6){};
        \node (n1) [blackV] at (0.3,-0.6){};
        \node (n1) [blackV] at (-0.3,-0.6){};

        \draw [line width=1.5pt] (-0.5,0.3) to (0.5,0.3);
        \draw [line width=1.5pt] (-0.5,-0.3) to (0.5,-0.3);
    \end{scope}
    \begin{scope}[xshift=4.5cm]
        \node (t1) [blackV] at (0,1){};
        \node (t2) [blackV] at (0,0){};
        \node (t3) [blackV] at (0,-1){};
        \node (s1) [blackV] at (1.25,1.5){};
        \node (s2) [blackV] at (1.25,0.5){};
        \node (s3) [blackV] at (1.25,-0.5){};
        \node (s4) [blackV] at (1.25,-1.5){};
        \node (l2) at (1.75,-2){$T_{1111}$};

        \foreach \i/\j in {t1/t2,t2/t3,s2/t2,t2/s3}
            \draw [arc] (\i) to (\j);
        \draw [arc, bend right=25] (t1) to (t3);
        \foreach \i in {1,2,3}{
            \draw [arc] (s1) to (t\i);
            \draw [arc] (t\i) to (s4);
        }
        \foreach \i in {2,3}{
            \draw [arc] (t1) to (s\i);
            \draw [arc] (s\i) to (t3);
        }
    \end{scope}
    \begin{scope}[xshift=7cm]
        \node (n1) [blackV] at (0,0.9){};
        \node (n1) [blackV] at (0,0.3){};
        \node (n1) [blackV] at (0,-0.3){};
        \node (n1) [blackV] at (0,-0.9){};

        \draw [line width=1.5pt] (-0.5,0.6) to (0.5,0.6);
        \draw [line width=1.5pt] (-0.5,0) to (0.5,0);
        \draw [line width=1.5pt] (-0.5,-0.6) to (0.5,-0.6);
    \end{scope}
    \begin{scope}[xshift=9cm]
        \node (t1) [blackV] at (0,1.5){};
        \node (t2) [blackV] at (0,0.5){};
        \node (t3) [blackV] at (0,-0.5){};
        \node (t4) [blackV] at (0,-1.5){};
        \node (s1) [blackV] at (1.25,2){};
        \node (s2) [blackV] at (1.25,0){};
        \node (s3) [blackV] at (1.25,-2){};
        \node (l3) at (1.75,-2.5){$T_{10101}$};

        \foreach \i/\j in {t1/t2,t2/t3,t3/t4,t1/s2,t2/s2,s2/t3,s2/t4}
            \draw [arc] (\i) to (\j);
        \foreach \i/\j in {t1/t3,t2/t4}
            \draw [arc, bend right=25] (\i) to (\j);
        \draw [arc, bend right=30] (t1) to (t4);
        \foreach \i in {1,2,3,4}{
            \draw [arc] (s1) to (t\i);
            \draw [arc] (t\i) to (s3);
        }
    \end{scope}
    \begin{scope}[xshift=11.5cm]
        \node (n1) [blackV] at (0,1.2){};
        \node (n1) [blackV] at (0,0){};
        \node (n1) [blackV] at (0,-1.2){};
        
        \draw [line width=1.5pt] (-0.5,0.9) to (0.5,0.9);
        \draw [line width=1.5pt] (-0.5,0.3) to (0.5,0.3);
        \draw [line width=1.5pt] (-0.5,-0.3) to (0.5,-0.3);
        \draw [line width=1.5pt] (-0.5,-0.9) to (0.5,-0.9);
    \end{scope}
\end{tikzpicture}
\caption{The family $\mathcal{F}_r$ of forbidden reflexive multipartite
tournaments.}
\label{fig:ref-mul-tour}
\end{figure}

Recall that reflexive strongly chordal tournaments are acyclic.    In view of
\Cref{lem:ref-comsplit}, the vertex set of any reflexive strongly chordal
multipartite tournament $D$ admits a partition $(T,S)$ where $T$ induces a
transitive tournament and $S$ is an independent set such that each of its
vertices is adjacent to every vertex in $T$.   Thus, for any $s \in S$, the
subdigraph of $D$ induced by $T \cup \{ s \}$ is again a transitive tournament,
so with respect to an acyclic ordering of $T$, we have that $s$ is an
out-neighbour of an initial segment and an in-neighbour of a terminal segment of
the acyclic ordering (see \Cref{fig:ref-mul-tour}).   So, every reflexive
strongly chordal tournament is completely determined by the number of vertices
in $T$ and the out-degree sequence of vertices in $S$.   From the previous
description of the structure of $D$, we also notice that vertices in $S$ can be
organized in $|T|+1$ ``levels'' (starting from the $0$-th level), where vertices
in the $i$-th level have exactly $i$ in-neighbours in $T$.   This description
works not only for reflexive strongly chordal tournaments, but for any
transitive orientation of a complete split graph; in \Cref{fig:ref-mul-tour} we
show examples of such oriented graphs, together with diagrams capturing all the
information needed to reconstruct the oriented graph: the number of lines in the
diagram corresponds to the number of vertices in $T$, and dots between the lines
represent the vertices in each of the $|T|+1$ different levels. Denote the level
of a vertex $x$ in $S$ by $\ell(x)$.

\begin{lemma}
\label{lem:split-ord}
    Let $D$ be a transitive orientation of a split graph with split partition
    $(T,S)$.
    \begin{enumerate}
        \item If a vertex in $T$ has at least two neighbours in $S$, then it
            cannot be placed before such neighbours in a $\Gamma$-free ordering.
        \item For any vertex $x$ in $S$, vertices in $N^-(x)$ appearing after
            $x$ in a $\Gamma$-free ordering of $D$ must appear in increasing
            order with respect to their levels, and vertices in $N^+(x)$
            appearing after $x$ in a $\Gamma$-free ordering of $D$ must appear
            in decreasing order with respect to their levels.
    \end{enumerate}
\end{lemma}

\begin{proof}
    The first item follows directly from \Cref{lem:ref-or-P3}.   The second item
    follows directly from \Cref{lem:simple-ord}.
\end{proof}

\begin{lemma}
\label{lem:ref-far-lev}
    Let $D$ be a transitive orientation of a complete split graph with split
    partition $(T,S)$, and let $x$ be a vertex in $S$.   If there are at least
    two vertices $y$ and $z$ in $S$ such that $|\ell(x) - \ell(y)|, |\ell(x) -
    \ell(z)| \ge 2$, then $x$ cannot be placed before $y$ and $z$ in a
    $\Gamma$-free ordering of $D$.
\end{lemma}
\begin{proof}
    Proceeding by contrapositive, we show that any ordering in which $x$ is
    placed before $x$ and $y$ is not $\Gamma$-free.   Suppose first that
    $\ell(x) < \ell(y)$ and $\ell(x) < \ell(z)$. By hypothesis, there are
    vertices $u$ and $v$ in $T$ which are out-neighbours of $x$ and
    in-neighbours of $y$ and $z$.   It follows from the first item of
    \Cref{lem:split-ord} that $x$ and $y$ or $x$ and $z$ must appear before $u$
    and $v$ in any $\Gamma$-free ordering of $D$.   But now, any ordering
    placing $x$ and $y$ or $z$ before $u$ and $v$ is not $\Gamma$-free, as
    otherwise the second item of \Cref{lem:split-ord} would imlpy that both $u$
    must appear before $v$ and $v$ must appear before $u$.

    Now suppose that $\ell(y) < \ell(x) < \ell(z)$.   By hypothesis, there are
    vertices $u$ and $v$ in $T$ which are out-neighbours of $y$ and
    in-neighbours of $x$, and verties $u'$ and $v'$ in $T$ which are
    out-neighbours of $x$ and in-neighbours of $z$.   Again, the first item of
    \Cref{lem:split-ord} implies that $x$ and $y$, or $x$ and $z$ must appear
    before $u$ and $v$ in any $\Gamma$-free ordering of $D$.   But in any case,
    the obtained order is not $\Gamma$-free, as otherwise the second item of
    \Cref{lem:split-ord} would create a contradiction with the order of $u$ and
    $v$, or the order of $u'$ and $v'$.

    The remaining cases, this is $\ell(z) < \ell(x) < \ell(y)$ and $\ell(y),
    \ell(z) < \ell(x)$, are similar, so their proofs are omitted.
\end{proof}

\begin{lemma}
\label{lem:ref-mul-mob}
    The three oriented graphs shown in \Cref{fig:ref-mul-tour} are multipartite
    tournament minimal obstructions for strong chordality.
\end{lemma}
\begin{proof}
    It follows from \Cref{lem:ref-or-P3,lem:ref-far-lev} that no vertex in
    $T_{202}$ and $T_{10101}$ can be the first vertex in a $\Gamma$-free
    ordering, and hence these oriented graphs are not strongly chordal. It also
    follows from the same lemmas that all vertices in the transitive tournament
    part, as well as vertices in levels $0$ and $3$ in the stable part of
    $T_{1111}$, cannot be placed first in a $\Gamma$-free ordering.
    
    Call $s_i$ the only vertex in the $i$-th level in the stable part, and let
    $(t_1, t_2, t_3)$ be a topological ordering of the transitive tournament
    part of $T_{1111}$ with $t_1$ a source.   We have already noticed that any
    $\Gamma$-free ordering of $T_{1111}$ must start with $s_1$ or $s_2$.   We
    only verify that $s_1$ cannot be the first vertex in a $\Gamma$-free
    ordering, the proof for $s_2$ is analogous.   If $s_1$ is first in a linear
    ordering of the vertices of $T_{1111}$, then we deduce from
    \Cref{lem:split-ord} that $t_2$ must be placed before $t_3$ and thus that
    $s_3$ cannot be placed before $t_2$ and $t_3$ in a $\Gamma$-free ordering of
    $T_{1111}$.   We obtain from \Cref{lem:ref-or-P3} that $s_0$ and $s_2$ must
    be placed before vertices in $T$ in a $\Gamma$-free ordering of $T_{1111}$.
    Hence, the second and third vertices in a $\Gamma$-free ordering of
    $T_{1111}$ must be $s_0$ and $s_2$ (in some order), nonetheless, an
    application of \Cref{lem:split-ord} creates a contradiction with the order
    in which $t_1$ and $t_2$ appear in a possible $\Gamma$-free ordering.
    Therefore, $T_{1111}$ is not strongly chordal.

    To prove minimality for $T_{1111}$, we only provide $\Gamma$-free orderings
    when we delete $t_1, t_2, s_0$ or $s_1$, the remaining cases are analogous.
    When we delete $t_1$ the ordering $(s_0, s_1, s_2, t_2, t_3, s_3)$ is a
    $\Gamma$-free ordering.   When we delete $t_2$ we use the ordering $(s_0,
    s_1, s_2, t_1, t_3, s_3)$.   When $s_0$ is missing, a $\Gamma$-free ordering
    is $(s_3, s_2, t_3, t_2, s_1, t_1)$.   Finally, when $s_1$ is missing the
    desired ordering is $(s_3, s_2, t_3, t_2, t_1, s_0)$.

    Verifications of minimality for $T_{202}$ and $T_{10101}$ are analogous.
\end{proof}

We are ready to prove the main result of this section.

\begin{theorem}
\label{thm:ref-mul-cha}
    Let $D$ be a reflexive multipartite tournament.   The following statements
    are equivalent.
    \begin{enumerate}
        \item $D$ is strongly chordal.
        \item $D$ is an $\mathcal{F}_r$-free transitive orientation of a
            complete split graph.
        \item $D$ is a transitive orientation of a complete split graph such
            that all vertices in $S$, except possible one, lie in two
            consecutive levels.
    \end{enumerate}
\end{theorem}
\begin{proof}
    It follows from \Cref{lem:ref-comsplit,lem:ref-mul-mob} that the first item
    implies the second one.   Suppose that $D$ is an $\mathcal{F}_r$-free
    transitive orientation of a complete split graph.   Being $T_{202}$-free,
    there are no non-consecutive levels with more than one vertex.   Having
    $T_{1111}$ as a forbidden induced subdigraph enforces $D$ to have at most
    three non-empty levels.   If there are three non-empty levels, then the
    absence of $T_{10101}$ as an induced subdigraph implies that there is a pair
    of consecutive levels in $D$.   Hence, all conditions of the third item are
    met.

    Finally, suppose that $D$ meets all conditions in the last item.  Suppose
    that $(t_1, \dots, t_n)$ is a topological ordering of the transitive
    tournament part of $D$, where $t_1$ is a source.    We first consider the
    case when there are at three non-empty levels, say levels $i, j$ and $k$,
    with $i+1 = j < k$.   Suppose that the vertices in the $i$-th and $j$-th
    levels are $s_{i1}, \dots, s_{ik_i}$, and $s_{j1}, \dots, s_{jk_j}$,
    respectively. By hypothesis there is only one vertex in level $k$, say
    $s_k$. Consider the ordering $(s_{i1}, \dots, s_{ik_i}, s_{j1}, \dots,
    s_{jk_j}, t_j, \dots, t_1, t_{j+1}, \dots, t_k, s_k, t_{k+1}, \dots, t_n)$.
    Observe that the adjacency matrix associated to the proposed ordering can be
    divided in 9 blocks: top left block is an identity matrix; top center is a
    constant $0$ matrix, except for the first $k_i$ rows of the first column;
    top right block is a constant $1$ matrix, except for the column
    corresponding to $s_k$; middle left block is a constant $1$ matrix, except
    for the first $k_i$ columns of the first row; middle center block is a lower
    triangular matrix having all entries on and below the main diagonal equal to
    $1$; middle right block is a constant $1$ matrix; bottom left and bottom
    center blocks are constant $0$ matrices; bottom right block is an upper
    triangular matrix having all entries on and above the main diagonal equal to
    $1$.   See \Cref{fig:ref-mul-cha} for an example of the adjacency matrix
    associated with the proposed ordering.   It is clear from the structure of
    the obtained adjacency matrix that the proposed ordering is $\Gamma$-free.
    
    The remaining cases are analogous to this one, or are obtained from it by
    deleting vertices (probably all) in some of the levels, so the restriction
    of the proposed ordering still works.
\end{proof}

\begin{figure}
    \footnotesize{
    \[
    \begin{pNiceArray}{cccc|cccc|ccc}[first-row,first-col]
     & s_{31} & s_{32} & s_{41} & s_{42} & t_4 & t_3 & t_2 & t_1 & t_5 & s_5  &t_6 \\ 
    s_{31} & 1 & 0 & 0 & 0 & 1 & 0 & 0 & 0 & 1 & 0 & 1 \\
    s_{32} & 0 & 1 & 0 & 0 & 1 & 0 & 0 & 0 & 1 & 0 & 1 \\
    s_{41} & 0 & 0 & 1 & 0 & 0 & 0 & 0 & 0 & 1 & 0 & 1 \\
    s_{42} & 0 & 0 & 0 & 1 & 0 & 0 & 0 & 0 & 1 & 0 & 1 \\
    \hline
    t_4    & 0 & 0 & 1 & 1 & 1 & 0 & 0 & 0 & 1 & 1 & 1 \\
    t_3    & 1 & 1 & 1 & 1 & 1 & 1 & 0 & 0 & 1 & 1 & 1 \\
    t_2    & 1 & 1 & 1 & 1 & 1 & 1 & 1 & 0 & 1 & 1 & 1 \\
    t_1    & 1 & 1 & 1 & 1 & 1 & 1 & 1 & 1 & 1 & 1 & 1 \\
    \hline
    t_5    & 0 & 0 & 0 & 0 & 0 & 0 & 0 & 0 & 1 & 1 & 1 \\
    s_5    & 0 & 0 & 0 & 0 & 0 & 0 & 0 & 0 & 0 & 1 & 1 \\
    t_6    & 0 & 0 & 0 & 0 & 0 & 0 & 0 & 0 & 0 & 0 & 1 \\
    \end{pNiceArray}
    \]
    }
    \caption{An example for the ordering used in the proof of
    \Cref{thm:ref-mul-cha}.}
    \label{fig:ref-mul-cha}
\end{figure}

We now use the nice structure of strongly chordal reflexive multipartite
tournaments to obtain a linear-time recognition algorithm.

\begin{corollary}
    There is an $O(|V|+|A|)$-time certifying algorithm to recognize strongly
    chordal reflexive multipartite tournaments.
\end{corollary}
\begin{proof}
    Given an input digraph $D$, we can test it to be reflexive, build its
    underlying graph, and test it to be a complete split graph in time
    $O(|V|+|A|)$.   If any of these tests fails, we can return either an
    irreflexive vertex, an induced $K_1 + K_2$ or an induced $C_4$ with the same
    time complexity.   If the tests succeed, we obtain the complete split
    partition of $D$.

    We can find a vertex of maximum out-degree in $O(|V|)$-time, and using DFS
    we can determine whether $D$ is acyclic in time $O(|V|+|A|)$, returning a
    directed cycle if it is not.   Finally, by inspecting only the out-degree of
    each vertex in the stable set of the partition, we can determine the level
    of each vertex, and verify whether one of $T_{202}, T_{1111}$ or $T_{10101}$
    occurs as an induced subdigraph of $D$, all this in time $O(|V|)$.   If the
    test fails, we can return one of such graphs as a no-certificate. Otherwise,
    we return the ordering given by the proof of \Cref{thm:ref-mul-cha} as a
    yes-certificate.
\end{proof}


\section{An example application}
\label{sec:application}

Domination in graphs is a well-known and widely studied subject.   Nonetheless,
its directed counterpart has received a lot less attention.   As one can notice
in the classic book of Haynes, Hedetniemi and Slater \cite{haynes1998}, the
study of domination in directed graphs at the end of the nineties was very
restricted.   This tendency has not changed a lot in the last 25 years, although
many variants of domination have now been considered in the literature at some
point; take as an example total and connected domination in digraphs
\cite{arumugamAJC39}.   Another domination related concept, the maximum number
of pairwise disjoint (either closed or open) in-neighbourhoods of $D$ was
recently considered in \cite{BresarJGT99}.   It is not hard to observe that a
single version of domination that behaves like the regular version of domination
or total domination, depending on the prescense or abscence of loops, can be
defined for digraphs with possible loops.  The same applies for the maximum
number of pairwise disjoint in-neighbourhoods (where loops are used to
distinguish between the open and closed neighbourhoods).

In \cite{hhhl} we explained how to unify the problems of domination and total
domination in graphs, by using loops to indicate whether a vertex can or cannot
dominate itself, and gave an algorithm to solve the unified problem on strongly
chordal graphs with possible loops. Here we observe that a small modification of
the algorithm actually solves the more general domination problem in strongly
chordal digraphs given with a strong ordering.

In the {\em digraph domination problem} for a digraph $G$ we seek to minimize
the number of vertices whose out-neighbourhoods cover all vertices, i.e., find a
smallest set of vertices $D$ such that each vertex in $G$ is dominated by some
vertex in $D$. A set $D$ with this property is called a {\em dominating set} of 
vertices. Note that here again, a vertex $v$ which has a loop can dominate
itself, while a loopless vertex can not. In the {\em disjoint in-neighbourhood
problem} we seek a largest set of vertices whose in-neighbourhoods are disjoint.
Since $k$ vertices with disjoint in-neighbourhoods need at least $k$ vertices to
dominate them, we have $|C| \leq |D|$ for any dominating set of vertices $D$ and
set of vertices $C$ with disjoint in-neighbourhoods, and if we find a dominating
set $D$ and a set $C$ of vertices with disjoint in-neighbourhoods such that
$|C|=|D|$, then $D$ is minimum and $C$ is maximum.

Thus assume we have a strongly chordal digraph $G$ presented with a strong
ordering $<$. We now explain how to compute such sets $C$ and $D$ iteratively.
As in \cite{hhhl} this is inspired by the algorithm of Martin Farber
\cite{domdom}, and motivated by the duality and complementary slackness of
linear programming. We use labels $N$ (indicating an unavailable vertex), and
$D$ and $C$ (for vertices forming the set $D$ and $C$). We assume that $G$ has
no vertices with in-degree zero, as such vertices cannot be dominated and the
problem has no solution.

\begin{itemize}
    \item Find, in the ordering $<$, the first vertex $x$ without the label $N$.
    \item Find, in the ordering $<$, the last in-neighbour $y$ of $x$.
    \item Label $x$ by $C$, label $y$ by $D$, and label all out-neighbours of
        $y$ by $N$.
\end{itemize}

A vertex will in general receive several labels, and each vertex will receive at
least one of the labels $C, N$. Moreover, every vertex will receive the label
$C$ and $D$ at most once. Specifically, when a vertex $x$ is labeled $C$, a
unique in-neighbour $y$ of $x$ is labeled $D$. At that point, all out-neighbours
of $y$, including $x$, receive the label $N$. Therefore $x$ will not receive the
label $C$ again. Similarly, $y$ will never receive another label $D$, since all
its out-neighbours are ineligible for the label $C$. We now have some vertices
labeled $C$ and the same number of vertices labeled $D$. Some vertices may have
both labels $C$ and $D$, but $|C|=|D|$. Moreover, the set $D$ is dominating, as
there are no vertices left without a label $C$ or $N$, and each vertex with
label $C$ or $N$ has a in-neighbour labeled $D$. Finally, we note that the
in-neighbourhoods of vertices labeled $C$ are disjoint. Otherwise some $x < x'$
both labeled $C$ have a common in-neighbour $z$. Suppose $y$ was the last
in-neighbour of $x$ when $x$ was labeled $C$: then we have $z \le y$. Since $x'$
is labeled $C$ later than $x$, it is not an in-neighbour of $y$ and hence $z$ is
not $y$. Therefore we have $x < x', z < y$, and $zx \in E(G), yx \in E(G)$, and
$zx' \in E(G), yx' \notin E(G)$, forming a $\Gamma$ submatrix, contradicting the
fact that $<$ is a strong ordering. Thus we have a general dominating set $D$
and a set of vertices $C$ with disjoint in-neighbourhoods, and $|C|=|D|$ whence
both are optimal.

This is a linear time algorithm solving the digraph domination problem (and the
disjoint in-neighbourhoods problem) if a strong ordering of the digraph is
given. Moreover our discussion proves the following minmax theorem.

\begin{theorem}
In a strongly chordal digraph the minimum number of vertices in a dominating
set equals the maximum number of vertices with disjoint in-neighbourhoods.
\end{theorem}

\section{Conclusions}
\label{sec:conclusions}
We have seen that strongly chordal digraphs can be recognized in polynomial time
amongst tournaments with possible loops, reflexive complete multipartite
tournaments, irreflexive complete bipartite tournaments, and irreflexive
complete multipartite tournaments minus one edge. We do not know if strongly
chordal digraphs can be recognized in polynomial time for a general digraph with
possible loops.

\begin{problem}
    Determine whether strongly chordal digraphs with possible loops can be
    recognized in polynomial time.
\end{problem}

We now mention one other natural class of digraphs with polynomial recognition
of strong chordality. Each bipartite graph $G$ defines a digraph $D_G$ by
orienting all edges from red to blue vertices; the adjacency matrix of $D_G$ is
clearly obtained from the bi-adjacency matrix of $G$ by adding rows and columns
of zeros. Thus independent permutations of rows and columns of $N(G)$ again
yield a symmetric ordering of $M(D_G)$. This means that $G$ is a chordal bigraph
if and only if $D_G$ is a strongly chordal digraph.

A {\em balanced digraph} is a digraph $D$ such that any cycle has the same
number of forward and backward arcs. By definition, a balanced digraph $D$ is
irreflexive, and it is easy to see that there is a vertex partition into parts
$V_i, i\in \{1, \dots, k\}$, such that each arc of $D$ starts in some $V_i$
and ends in $V_{i+1}$. The adjacency matrix of a balanced digraph can be
symmetrically permuted into consecutive blocks corresponding to the parts $V_i$.
In such a form, a symmetric permutation of the matrix corresponds to independent
permutations of rows and columns in each submatrix $M_i$ with rows in block
$V_i$ and columns in block $V_{i+1}$. Moreover, it is easy to see that each
$\Gamma$ submatrix of $M$ must lie in some $M_i$. Note that when $k=2$, i.e.,
when there are only two parts, $V_1, V_2$, a balanced digraph is some $D_G$ for
a bipartite graph $G$. For a general balanced digraph, denote by $G_i$ the
underlying bipartite subgraph of $D$ with parts $V_i, V_{i+1}$.

\begin{theorem}
A balanced digraph $D$ is strongly chordal if and only if each $G_i$ is a
chordal bigraph.
\end{theorem}

We can translate this result to a forbidden subgraph characterization. A {\em
fence} is an oriented even cycle of length greater than four, without a directed
path of length two, see Figure \ref{fig:fff}.

\begin{corollary}
A balanced digraph $D$ is strongly chordal if and only if it does not contain a
fence as an induced subgraph.
\end{corollary}

\begin{figure}[htb!]
\centering
\begin{tikzpicture}
    \begin{scope}
        \node (1u) [vertex] at (-1.5,1){};
        \node (2u) [vertex] at (0,1){};
        \node (pu) at (1.3,1){$\cdots$};
        \node (3u) [vertex] at (2.5,1){};
    \end{scope}
    \begin{scope}[xshift=0.75cm]
        \node (1d) [vertex] at (-1.5,-0.25){};
        \node (2d) [vertex] at (0,-0.25){};
        \node (pd) at (1.2,-0.25){$\cdots$};
        \node (3d) [vertex] at (2.5,-0.25){};
    \end{scope}
    \foreach \u/\v in {1d/1u,1d/2u,2d/2u,3d/1u,3d/3u}
        \draw [arc] (\u) to (\v);
\end{tikzpicture}
\caption{A fence.} \label{fig:fff}
\end{figure}

\begin{corollary}
Each oriented tree is strongly chordal.
\end{corollary}

\end{document}